\documentclass[manuscript,printscheme,reqno]{aomart}
\usepackage{mystyle}
\usepackage{mathtools}
\usepackage{stmaryrd}
\usepackage{amssymb}
\usepackage{letltxmacro}

\newtheorem{thm}{Theorem}[section]
\newtheorem{cor}[thm]{Corollary}
\newtheorem{lem}[thm]{Lemma}
\newtheorem{prop}[thm]{Proposition}
\theoremstyle{definition}
\newtheorem{defn}[thm]{Definition}
\newtheorem{rem}[thm]{Remark}
\LetLtxMacro\origproof\proof
\LetLtxMacro\endorigproof\endproof
\makeatletter
\renewenvironment{proof}[1][\proofname]{\origproof[#1]\mbox{}\par\nobreak\@afterindentfalse\@afterheading\ignorespaces}{\endorigproof}
\makeatother

\DeclarePairedDelimiter{\floor}{\lfloor}{\rfloor}
\DeclarePairedDelimiter{\abs}{\lvert}{\rvert}
\DeclareMathOperator{\lcm}{lcm}
\newcommand{\N}{\mathbb{N}}
\newcommand{\Z}{\mathbb{Z}}
\newcommand{\R}{\mathbb{R}}
\newcommand{\bbinom}[3][b]{\genfrac{\llbracket}{\rrbracket}{0pt}{}{#2}{#3}_{\!#1}}

\allowdisplaybreaks

\title[Denumerant II]{Arithmetic Upper and Lower Bounds for the Denumerant Function II}
\author{Gerardo Ballesio}
\email{g.ballesio@gmail.com}
\givenname{Gerardo}
\surname{Ballesio}
\copyrightnote{}
\subject{primary}{matsc2020}{11P81}
\subject{secondary}{matsc2020}{05A17}
\subject{secondary}{matsc2020}{11B34}

\begin{document}

\begin{abstract}
	We apply the elementary arithmetic argument of the first paper of this series
	(arXiv:2204.13689) to binomial coefficients, and obtain upper and lower bounds for Sylvester's
	denumerant function for arbitrary positive integers $a_1,\dotsc,a_k$ that sharpen those of
	Agnarsson. The lower bound is a sum of binomial coefficients; it holds along the multiples of
	$\gcd(a_1,\dotsc,a_k)$ without assuming that a solution exists, and it is greater than the
	polynomial lower bound of the first paper. In the upper bound the shift of Agnarsson's
	binomial coefficient is reduced.
\end{abstract}
\maketitle
\tableofcontents

\section*{Notation}\label{s:notation}
Throughout the paper we use the following notation.
\begin{itemize}
	\item $\N=\{0,1,2,\dotsc\}$ and, for $t\in\N$, $\N_{\geq t}=\{t,t+1,t+2,\dotsc\}$; $\Z$ and $\R$
	      are the integers and the reals, and $\R_{\geq0}$ the non-negative reals.
	\item For $x,y\in\Z$ with $x\neq0$, $x\mid y$ means that $y=cx$ for some $c\in\Z$.
	\item $\gcd$ and $\lcm$ are the greatest common divisor and the least common multiple of
	      positive integers, and $\floor{x}$ is the floor of $x\in\R$.
	\item For $m\in\Z$ and $j\in\N$, $\binom{m}{j}\coloneqq\frac{m(m-1)\dotsm(m-j+1)}{j!}$.
	\item For $m\in\N$ and $i\in\Z$, $e_i(X_1,\dotsc,X_m)$ is the elementary symmetric polynomial of
	      degree $i$, with $e_0=1$, and $e_i=0$ for $i<0$ or $i>m$.
	\item For functions $f,g$ on $\N$ with $g\geq0$, $f(n)=O\bigl(g(n)\bigr)$ means that there are
	      a real number $C>0$ and $n_0\in\N$ with $\abs{f(n)}\leq C\,g(n)$ for every $n\geq n_0$.
	\item Empty sums are $0$ and empty products are $1$.
	\item $a_1,a_2,a_3,\dotsc$ is a fixed sequence of positive integers. Only finitely many of
	      them are involved in each statement.
\end{itemize}

\section{Introduction}\label{s:intro}
In the nineteenth century, while studying the partition function,
J.~J.~Sylvester~\cite{sylvester:partition} and A.~Cayley~\cite{cayley} introduced
the \emph{denumerant}: the number of ways of writing a non-negative integer $n$ as a
non-negative integer combination of given positive integers $a_1,\dotsc,a_k$.

\begin{defn}\label{dfn:denumerant}
	For $k\in\N_{\geq1}$ and $n\in\N$, the \emph{denumerant} $D_k(n)$ is the number of
	solutions $(x_1,\dotsc,x_k)\in\N^k$ of $a_1x_1+\dotsb+a_kx_k=n$. Moreover, we set
	$d_k\coloneqq\gcd(a_1,\dotsc,a_k)$; in particular $d_1=a_1$.
\end{defn}

Since $d_k$ divides each of $a_1,\dotsc,a_k$, it divides every value of
$a_1x_1+\dotsb+a_kx_k$; so $D_k(n)=0$ unless $n$ is a multiple of $d_k$. This paper is a sequel
to~\cite{ballesio-arxiv}, which we call the first paper, and it gives new explicit upper and
lower bounds for $D_k(n)$.

\subsection{Known results}
We recall the two results that this paper starts from: the bounds of Agnarsson, in his
notation, and those of the first paper.

In 2002 G.~Agnarsson~\cite{agnarsson} proved, by elementary methods, that for $k\geq2$
\begin{equation}\label{eq:intro-agnarsson}
	\frac{d_k}{a_1\dotsm a_k}\binom{n-B_k}{k-1}
	\;\leq\; D_k(n) \;\leq\;
	\frac{d_k}{a_1\dotsm a_k}\binom{n+A_k}{k-1},
\end{equation}
where $A_1=B_1=0$ and, for $k\geq2$,
\begin{equation}\label{eq:intro-AB}
	A_k=A_{k-1}+\frac{d_{k-1}}{d_k}\,a_k,
	\qquad
	B_k=B_{k-1}+\biggl(\frac{d_{k-1}}{d_k}-1\biggr)a_k-1.
\end{equation}
In~\cite{agnarsson} the upper bound is proved for every $n\in\N$ (Theorem~3.3), and the
lower bound for $n\geq B_k+k-1$ under the assumption $D_k(n)>0$ (Theorem~3.4). Each of the two
bounds is a single binomial coefficient of degree $k-1$ in $n$.

In the first paper the bounds are powers of $n$ instead of binomial coefficients. They are
written in terms of the following integers.

\begin{defn}\label{dfn:constants}
	For $k\in\N_{\geq2}$ we set
	\begin{equation*}
		b_k\coloneqq\lcm(d_{k-1},a_k),
		\qquad
		\Delta_k\coloneqq\sum_{i=2}^{k}(b_i-a_i),
		\qquad
		\nabla_k\coloneqq b_2+\frac12\sum_{i=3}^{k}b_i .
	\end{equation*}
\end{defn}

The main result of the first paper (Theorem~1.3 there) is that, for $k\geq2$ and every multiple
$n$ of $d_k$,
\begin{equation}\label{eq:intro-first}
	\frac{d_k\,(n-\Delta_k)^{k-1}}{(k-1)!\,a_1\dotsm a_k}
	\;<\; D_k(n) \;\leq\;
	\frac{d_k\,(n+\nabla_k)^{k-1}}{(k-1)!\,a_1\dotsm a_k},
\end{equation}
where the lower bound holds for $n\geq\Delta_k$. For $a_1=1$ these are the bounds of G.~Blom and
C.~Fr\"oberg~\cite{blom-froberg}.

\subsection{Main results}
Here we run the argument of the first paper with binomial coefficients in place of powers; this
sharpens both bounds~\eqref{eq:intro-agnarsson} of Agnarsson. The new bounds are written with
the following symbol and integers.

\begin{defn}\label{dfn:coefficients}
	Let $c=(c_3,c_4,c_5,\dotsc)$ be a sequence of real numbers. For $m\in\N$ and $i\in\Z$ we set
	\begin{equation*}
		\bbinom[c]{m}{i}\coloneqq e_i\Bigl(\frac{c_3}{2},\dotsc,\frac{c_{m+2}}{2}\Bigr).
	\end{equation*}
	We denote by $b-1$ the sequence $c_l=b_l-1$, so that
	\begin{equation}\label{eq:coeff-b-1}
		\bbinom[b-1]{m}{i}=e_i\Bigl(\frac{b_3-1}{2},\dotsc,\frac{b_{m+2}-1}{2}\Bigr).
	\end{equation}
	Moreover, for $k\in\N_{\geq2}$ we set
	\begin{equation}\label{eq:M-def}
		M_k\coloneqq b_2+\sum_{i=3}^{k}\Bigl(\floor*{\frac{b_i}{2}}+1\Bigr).
	\end{equation}
\end{defn}

The symbol for the sequence $b-1$ gives the coefficients of our lower bound, and $M_k$ is the
shift of $n$ in our upper bound, in place of $A_k$. Our main result is the following.

\begin{thm}\label{thm:main}
	Let $k\in\N_{\geq2}$ and let $n\in\N$ be a multiple of $d_k$. Then
	\begin{equation}\label{eq:intro-main}
		\frac{d_k}{a_1\dotsm a_k}\sum_{i=0}^{k-2}
		\bbinom[b-1]{k-2}{i}\binom{n-B_k-i}{k-1-i}
		\;\leq\; D_k(n) \;\leq\;
		\frac{d_k}{a_1\dotsm a_k}\binom{n+M_k}{k-1},
	\end{equation}
	where the lower bound holds for $n\geq B_k+k-1$ and the upper bound for every such $n$.
\end{thm}

The term $i=0$ of the lower bound is the lower bound of Agnarsson, and the other terms are
non-negative; moreover the assumption $D_k(n)>0$ is not needed. In the upper bound
$M_k\leq A_k$.

\subsection{Organization of the paper}
Section~\ref{s:prelim} collects the tools common to both bounds: the recursion for the
denumerant, Agnarsson's constants in our notation, binomial coefficients truncated at negative
arguments, and their sums along an arithmetic progression. Sections~\ref{s:lower}
and~\ref{s:upper} are devoted to the lower and to the upper bound of Theorem~\ref{thm:main}:
each of them contains the lemmas needed only for its bound, the theorem, and the comparison
with the corresponding bound of Agnarsson; the upper bound is proved in a sharper form, with
terms of lower degree when some of $b_3,\dotsc,b_k$ are even.
Section~\ref{s:compare} compares the new bounds with those of the first paper, proving that the
new lower bound is greater than the sharper, polynomial lower bound given there, and with bounds
in terms of the least common multiple of $a_1,\dotsc,a_k$; it also compares the distance between
the upper and the lower bound in the first paper and in this one. Finally,
Section~\ref{s:final} contains two remarks, as in the first paper: a refinement of both bounds
at the last step of the induction, and bounds for the number of solutions of
$a_1x_1+\dotsb+a_kx_k\leq n$.

\section{Preliminaries}\label{s:prelim}

\subsection{The integers \texorpdfstring{$d_k$, $b_k$ and $\Delta_k$}{dk, bk and Delta k}}
We recall from the first paper the properties of the integers of
Definitions~\ref{dfn:denumerant} and~\ref{dfn:constants} that we need.

\begin{lem}[parts of Lemma~2.1 of~\cite{ballesio-arxiv}]\label{lem:basic}
	For every $k\geq2$:
	\begin{enumerate}
		\item[\textup{(i)}] $d_k=\gcd(d_{k-1},a_k)$; hence $d_k\mid d_{k-1}$, and
		      $b_k=\frac{d_{k-1}}{d_k}\,a_k$; in particular $b_k$ is a multiple of $a_k$ and
		      $b_k\geq a_k\geq1$;
		\item[\textup{(ii)}] $b_2\dotsm b_k=\frac{a_1\dotsm a_k}{d_k}$;
		\item[\textup{(iii)}] $\Delta_k$ is a non-negative integer, $\Delta_2=b_2-a_2$, and
		      $\Delta_{k+1}=\Delta_k+b_{k+1}-a_{k+1}$;
		\item[\textup{(iv)}] $D_k(n)=0$ for every $n\in\N$ that is not a multiple of $d_k$.
	\end{enumerate}
\end{lem}
The numbers $b_k$ and $\Delta_k$ depend on the order of $a_1,\dotsc,a_k$, while $D_k(n)$ does
not: every bound below may be applied to the most favourable ordering.

\subsection{The recursion}
Passing from $k$ to $k+1$, the solutions of $a_1x_1+\dotsb+a_{k+1}x_{k+1}=n$ are classified
according to the value $t=a_{k+1}x_{k+1}$. Only the values of $t$ for which $n-t$ is a multiple
of $d_k$ contribute, and the next lemma shows that they form an arithmetic progression with
difference $b_{k+1}$; we refer to~\cite{ballesio-arxiv} for the proof.

\begin{lem}[Lemma~2.2 of~\cite{ballesio-arxiv}]\label{lem:recursion}
	Let $k\in\N_{\geq2}$ and let $n\in\N$ be a multiple of $d_{k+1}$. Then there is exactly one
	integer $s$ such that
	\begin{equation*}
		0\leq s<b_{k+1},\qquad a_{k+1}\mid s,\qquad d_k\mid n-s ;
	\end{equation*}
	moreover $s\leq b_{k+1}-a_{k+1}$, and
	\begin{equation*}
		D_{k+1}(n)=\sum_{q=0}^{\floor*{\frac{n-s}{b_{k+1}}}}D_k\bigl(n-s-qb_{k+1}\bigr),
	\end{equation*}
	where every argument $n-s-qb_{k+1}$ is a non-negative multiple of $d_k$.
\end{lem}

\subsection{The constants of Agnarsson in our notation}
In the first paper (Section~5.2 there) the constants $A_k$ and $B_k$ of~\eqref{eq:intro-AB} are
written in our notation:
\begin{equation}\label{eq:AB-ours}
	A_k=\sum_{i=2}^{k}b_i,
	\qquad
	B_k=\Delta_k-(k-1).
\end{equation}
In particular, $n\geq B_k+k-1$ if and only if $n\geq\Delta_k$. Moreover
$\frac{d_k}{a_1\dotsm a_k}=\frac{1}{b_2\dotsm b_k}$ by Lemma~\ref{lem:basic}~(ii). Hence the bounds~\eqref{eq:intro-agnarsson} read
\begin{equation}\label{eq:agnarsson-ours}
	\frac{1}{b_2\dotsm b_k}\binom{n-B_k}{k-1}
	\;\leq\; D_k(n) \;\leq\;
	\frac{1}{b_2\dotsm b_k}\binom{n+A_k}{k-1}.
\end{equation}
Throughout the paper the coefficient is written as $\frac{1}{b_2\dotsm b_k}$.

\subsection{Truncated binomial coefficients}
Our bounds for $D_k(n)$ are binomial coefficients $\binom{t+r}{r}$, or combinations of them, where
$t$ depends on $n$, for instance $t=n-\Delta_k$. To state them for every $n$, and to sum them along
arithmetic progressions that reach negative values of $t$, we set them equal to $0$ for $t<0$.

\begin{defn}\label{dfn:trunc}
	For $r\in\N$ and $t\in\Z$ let
	\begin{equation*}
		f_r(t)\coloneqq
		\begin{cases}
			\displaystyle\binom{t+r}{r}=\frac{(t+1)(t+2)\dotsm(t+r)}{r!}, & t\geq0,\\[2ex]
			0, & t<0.
		\end{cases}
	\end{equation*}
\end{defn}
For instance, for $t\geq0$ we have $f_0(t)=1$, $f_1(t)=t+1$ and $f_2(t)=\frac{(t+1)(t+2)}{2}$.
The function $f_r$ is itself a denumerant: for $t\geq0$, $\binom{t+r}{r}$ is the number of solutions
of $x_1+\dotsb+x_{r+1}=t$ in $\N^{r+1}$, so $f_r(t)=D_{r+1}(t)$ for the sequence $a_1=a_2=\dotsb=1$.
Our bounds are thus denumerants of the sequence $1,1,\dotsc$, evaluated at shifted arguments.

\begin{lem}\label{lem:trunc}
	Let $r\in\N_{\geq1}$. Then:
	\begin{enumerate}
		\item[\textup{(i)}] $f_r(t)-f_r(t-1)=f_{r-1}(t)$ for every $t\in\Z$;
		\item[\textup{(ii)}] $f_{r-1}$ and $f_r$ are non-negative and non-decreasing on $\Z$;
		\item[\textup{(iii)}] for every $T\in\N$ and every integer $c\leq0$,
		      $\displaystyle\sum_{t=c}^{T}f_{r-1}(t)=f_r(T)$.
	\end{enumerate}
\end{lem}
\begin{proof}
	\begin{enumerate}
		\item[\textup{(i)}] If $t\geq1$, then $t$ and $t-1$ are both $\geq0$, and the claim is Pascal's
		      rule $\binom{t+r}{r}-\binom{t+r-1}{r}=\binom{t+r-1}{r-1}$. If $t=0$, then
		      $f_r(0)=\binom{r}{r}=1$, $f_r(-1)=0$ and $f_{r-1}(0)=\binom{r-1}{r-1}=1$, so both sides
		      equal $1$. If $t<0$, then $t-1<0$ too, and both sides equal $0$.
		\item[\textup{(ii)}] For $t\geq0$, $f_r(t)$ is a product of positive numbers divided by $r!$, hence
		      positive; for $t<0$ it is $0$. So $f_r\geq0$, and likewise $f_{r-1}\geq0$. By~(i),
		      $f_r(t)-f_r(t-1)=f_{r-1}(t)\geq0$ for every $t$, so $f_r$ is non-decreasing. If $r\geq2$,
		      the same argument with $r-1$ in place of $r$ shows that $f_{r-1}$ is non-decreasing; if
		      $r=1$, then $f_{r-1}=f_0$ is $0$ for $t<0$ and $1$ for $t\geq0$, which is non-decreasing.
		\item[\textup{(iii)}] By~(i), each term of the sum is a difference of two consecutive values of
		      $f_r$:
		      \begin{equation*}
		      	\sum_{t=c}^{T}f_{r-1}(t)
		      	=\bigl(f_r(c)-f_r(c-1)\bigr)+\bigl(f_r(c+1)-f_r(c)\bigr)+\dotsb+\bigl(f_r(T)-f_r(T-1)\bigr).
		      \end{equation*}
		      All the intermediate values cancel, and only $f_r(T)-f_r(c-1)$ remains. Since $c\leq0$,
		      we have $c-1<0$, hence $f_r(c-1)=0$.\qedhere
	\end{enumerate}
\end{proof}

For $c=0$, part~(iii) is the parallel summation identity $\sum_{t=0}^{T}\binom{r-1+t}{t}=\binom{r-1+T+1}{T}$ of~\cite[(5.9)]{concrete-math};
thanks to the truncation it also holds when the sum starts at a negative integer.

The next lemma collects two properties of the increments of $f_r$, which are used in the estimates
of Sections~\ref{s:lower-sums} and~\ref{s:upper-sums}: the first one for both estimates, the second
one (convexity) for the upper estimate.

\begin{lem}\label{lem:increments}
	Let $r\in\N_{\geq1}$ and $t\in\Z$. Then:
	\begin{enumerate}
		\item[\textup{(i)}] for every $u\in\N$, $\displaystyle f_r(t)-f_r(t-u)=\sum_{v=0}^{u-1}f_{r-1}(t-v)$;
		\item[\textup{(ii)}] for every $j\in\N$, $f_r(t-j)+f_r(t+j)\geq2f_r(t)$.
	\end{enumerate}
\end{lem}
In words: (i) the increment of $f_r$ over $u$ consecutive steps is the sum of the values of $f_{r-1}$
at the upper ends of the steps; (ii) $f_r$ is \emph{convex}: its value at $t$ is at most the average of
its values at two points symmetric with respect to $t$. For instance, for $r=2$ and $t=j=1$, (ii)
reads $f_2(0)+f_2(2)=1+6\geq2\cdot3=2f_2(1)$.
\begin{proof}
	\begin{enumerate}
		\item[\textup{(i)}] For $u=0$ both sides are $0$. For $u\geq1$ we go down from $t$ to $t-u$ one
		      step at a time, and use Lemma~\ref{lem:trunc}~(i) at each step:
		      \begin{align*}
		      	f_r(t)-f_r(t-u)
		      	&=\bigl(f_r(t)-f_r(t-1)\bigr)+\bigl(f_r(t-1)-f_r(t-2)\bigr)+\dotsb+\bigl(f_r(t-u+1)-f_r(t-u)\bigr)\\
		      	&=f_{r-1}(t)+f_{r-1}(t-1)+\dotsb+f_{r-1}(t-u+1)
		      	=\sum_{v=0}^{u-1}f_{r-1}(t-v).
		      \end{align*}
		\item[\textup{(ii)}] By~(i), applied once at the point $t+j$ and once at the point $t$, both with
		      $u=j$,
		      \begin{equation*}
		      	f_r(t+j)-f_r(t)=\sum_{v=0}^{j-1}f_{r-1}(t+j-v),
		      	\qquad
		      	f_r(t)-f_r(t-j)=\sum_{v=0}^{j-1}f_{r-1}(t-v).
		      \end{equation*}
		      For each $v$ we have $t+j-v\geq t-v$, and $f_{r-1}$ is non-decreasing
		      (Lemma~\ref{lem:trunc}~(ii)); so each term of the first sum is at least the corresponding
		      term of the second, and $f_r(t+j)-f_r(t)\geq f_r(t)-f_r(t-j)$. Moving $f_r(t)$ and
		      $f_r(t-j)$ to the other side gives~(ii).\qedhere
	\end{enumerate}
\end{proof}

\subsection{The coefficients}
Our lower bound is a combination of the functions $f_{k-1-i}$ with coefficients that are elementary
symmetric polynomials. At each step of the induction two such combinations are added, and their
coefficients combine by a rule of Pascal type (Lemma~\ref{lem:pascal-sums}).

These coefficients are the symbols of Definition~\ref{dfn:coefficients}. We shall mainly use
two sequences: $c=b$, for which we omit the index and write $\bbinom{m}{i}$ as in Section~2.4 of
the first paper, and $c=b-1$ (see~\eqref{eq:coeff-b-1}). The second one gives the
coefficients of our lower bound; the first one is needed only in Section~\ref{s:compare}. A
third sequence is introduced in Section~\ref{s:upper}. For instance,
\begin{equation*}
	\bbinom{m}{1}=\sum_{l=3}^{m+2}\frac{b_l}{2},
	\qquad
	\bbinom[b-1]{m}{1}=\sum_{l=3}^{m+2}\frac{b_l-1}{2}.
\end{equation*}

\begin{lem}[Lemma~2.8~(iv) of~\cite{ballesio-arxiv}]\label{lem:sym-poly}
	Let $c=(c_3,c_4,\dotsc)$ be a sequence of real numbers, $m\in\N$ and $i\in\Z$. Then:
	\begin{enumerate}
		\item[\textup{(i)}] $\bbinom[c]{m}{0}=1$, and $\bbinom[c]{m}{i}=0$ for $i<0$ or $i>m$;
		\item[\textup{(ii)}] if $c_l\geq0$ for every $l$, then $\bbinom[c]{m}{i}\geq0$;
		\item[\textup{(iii)}] if $m\geq1$, then
		      $\bbinom[c]{m}{i}=\bbinom[c]{m-1}{i}+\frac{c_{m+2}}{2}\bbinom[c]{m-1}{i-1}$.
	\end{enumerate}
	In particular, since $b_l\geq1$ by Lemma~\ref{lem:basic}~(i), the numbers $b_l$ and $b_l-1$ are
	non-negative, and \textup{(ii)} applies to both $c=b$ and $c=b-1$.
\end{lem}

In the inductive steps of the proofs of Theorems~\ref{thm:main-lower}, \ref{thm:main-upper}
and~\ref{thm:compare-lower}, a sum with coefficients $\bbinom[c]{m-1}{i}$ and a second, shifted sum with the same coefficients are
added and combined into a single sum with coefficients $\bbinom[c]{m}{i}$. This is the content of the
next lemma, which is Lemma~\ref{lem:sym-poly}~(iii) written for sums.

\begin{lem}\label{lem:pascal-sums}
	Let $c=(c_3,c_4,\dotsc)$ be a sequence of real numbers, let $m\in\N_{\geq1}$, and let
	$y_0,y_1,\dotsc,y_m$ be real numbers. Then
	\begin{equation*}
		\sum_{i=0}^{m-1}\bbinom[c]{m-1}{i}y_i+\frac{c_{m+2}}{2}\sum_{i=0}^{m-1}\bbinom[c]{m-1}{i}y_{i+1}
		=\sum_{i=0}^{m}\bbinom[c]{m}{i}y_i .
	\end{equation*}
\end{lem}
For instance, for $m=1$ both sides equal $y_0+\frac{c_3}{2}\,y_1$, since
$\bbinom[c]{0}{0}=\bbinom[c]{1}{0}=1$ and $\bbinom[c]{1}{1}=\frac{c_3}{2}$.
\begin{proof}
	We start from the right-hand side. By Lemma~\ref{lem:sym-poly}~(iii), for every $0\leq i\leq m$,
	$\bbinom[c]{m}{i}=\bbinom[c]{m-1}{i}+\frac{c_{m+2}}{2}\bbinom[c]{m-1}{i-1}$, so
	\begin{equation*}
		\sum_{i=0}^{m}\bbinom[c]{m}{i}y_i
		=\sum_{i=0}^{m}\bbinom[c]{m-1}{i}y_i+\frac{c_{m+2}}{2}\sum_{i=0}^{m}\bbinom[c]{m-1}{i-1}y_i .
	\end{equation*}
	\begin{itemize}
		\item In the first sum, the term $i=m$ vanishes, since $\bbinom[c]{m-1}{m}=0$
		      (Lemma~\ref{lem:sym-poly}~(i)); so it is $\sum_{i=0}^{m-1}\bbinom[c]{m-1}{i}y_i$.
		\item In the second sum, the term $i=0$ vanishes, since $\bbinom[c]{m-1}{-1}=0$. Putting
		      $i=j+1$, with $0\leq j\leq m-1$, it becomes $\sum_{j=0}^{m-1}\bbinom[c]{m-1}{j}y_{j+1}$.
	\end{itemize}
	Renaming $j$ as $i$ gives the left-hand side.
\end{proof}

\subsection{Sums along an arithmetic progression}\label{s:sums}
For $r,w\in\N_{\geq1}$ and $T\in\Z$ we consider the sum of $f_r$ along the arithmetic progression
that starts at $T$ and goes down with difference $w$:
\begin{equation}\label{eq:S-def}
	S_{r,w}(T)\coloneqq\sum_{q\geq0}f_r(T-qw).
\end{equation}
Only finitely many terms are different from $0$, since $f_r(T-qw)=0$ as soon as $T-qw<0$. Thus
$S_{r,w}(T)=\sum_{q=0}^{\floor{T/w}}f_r(T-qw)$ for $T\in\N$, and $S_{r,w}(T)=0$ for $T<0$.
Separating the term $q=0$, and then putting $q=p+1$ in the remaining sum, so that
$T-qw=(T-w)-pw$, for every $T\in\Z$ we get
\begin{equation}\label{eq:S-rec}
	S_{r,w}(T)
	=f_r(T)+\sum_{q\geq1}f_r(T-qw)
	=f_r(T)+\sum_{p\geq0}f_r\bigl((T-w)-pw\bigr)
	=f_r(T)+S_{r,w}(T-w).
\end{equation}
The difference of the progression is denoted by $w$, to avoid confusion with the numbers
$b_2,b_3,\dotsc$; in the applications, $w=b_{k+1}$. The sum $S_{r,w}(T)$ is estimated from below
in Section~\ref{s:lower-sums} and from above in Section~\ref{s:upper-sums}.

\section{The lower bound}\label{s:lower}
Agnarsson's lower bound in~\eqref{eq:agnarsson-ours} is a single binomial coefficient.
Theorem~\ref{thm:main-lower} refines it by adding binomial coefficients of lower degree. Its
proof needs a lower estimate for the sums~\eqref{eq:S-def}, which we prove first.

\subsection{A lower estimate for the sums along a progression}\label{s:lower-sums}
Let $t\in\Z$. Since $f_r$ is non-decreasing (Lemma~\ref{lem:trunc}~(ii)), $f_r(t)$ is the largest
value of $f_r$ on the block $t,t-1,\dotsc,t-w+1$ below $t$, so
\begin{equation}\label{eq:block-agnarsson}
	w\,f_r(t)\geq\sum_{u=0}^{w-1}f_r(t-u).
\end{equation}
This is the estimate used by Agnarsson for each block in the proof of Lemma~3.2 of~\cite{agnarsson}.
The next lemma improves it.

\begin{lem}\label{lem:block-lower}
	Let $r,w\in\N_{\geq1}$ and $t\in\Z$. Then
	\begin{equation*}
		w\,f_r(t)\geq\sum_{u=0}^{w-1}f_r(t-u)+\frac{w-1}{2}\sum_{u=0}^{w-1}f_{r-1}(t-u).
	\end{equation*}
\end{lem}
\begin{proof}
	We prove the statement by induction on $w$.

	\emph{Base case $w=1$.} For $w=1$ the first sum is $f_r(t)$ and the factor $\frac{w-1}{2}$ is
	$0$, so the statement reads $f_r(t)\geq f_r(t)$.

	\emph{Inductive step.} Let $w\geq1$, and assume that
	\begin{equation}\label{eq:block-hyp}
		w\,f_r(t)\geq\sum_{u=0}^{w-1}f_r(t-u)+\frac{w-1}{2}\sum_{u=0}^{w-1}f_{r-1}(t-u).
	\end{equation}
	We have to prove that
	\begin{equation}\label{eq:block-goal}
		(w+1)\,f_r(t)\geq\sum_{u=0}^{w}f_r(t-u)+\frac{w}{2}\sum_{u=0}^{w}f_{r-1}(t-u).
	\end{equation}
	By Lemma~\ref{lem:increments}~(i) with $u=w$,
	$f_r(t)=f_r(t-w)+\sum_{u=0}^{w-1}f_{r-1}(t-u)$. Each of the $w$ terms of this sum is at least
	$f_{r-1}(t-w)$, since $f_{r-1}$ is non-decreasing (Lemma~\ref{lem:trunc}~(ii)); so half of the
	sum is at least $\frac{w}{2}\,f_{r-1}(t-w)$. Hence
	\begin{equation}\label{eq:block-one}
		f_r(t)
		=f_r(t-w)+\frac12\sum_{u=0}^{w-1}f_{r-1}(t-u)+\frac12\sum_{u=0}^{w-1}f_{r-1}(t-u)
		\geq f_r(t-w)+\frac12\sum_{u=0}^{w-1}f_{r-1}(t-u)+\frac{w}{2}\,f_{r-1}(t-w).
	\end{equation}
	Adding~\eqref{eq:block-hyp} and~\eqref{eq:block-one}, and then collecting the terms,
	\begin{align*}
		(w+1)\,f_r(t)
		&\geq\sum_{u=0}^{w-1}f_r(t-u)+f_r(t-w)
		+\Bigl(\frac{w-1}{2}+\frac12\Bigr)\sum_{u=0}^{w-1}f_{r-1}(t-u)+\frac{w}{2}\,f_{r-1}(t-w)\\
		&=\sum_{u=0}^{w}f_r(t-u)+\frac{w}{2}\sum_{u=0}^{w}f_{r-1}(t-u),
	\end{align*}
	which is~\eqref{eq:block-goal}.
\end{proof}

\begin{lem}[cf.\ Lemma~3.2 of~\cite{agnarsson}]\label{lem:sum-lower}
	Let $r,w\in\N_{\geq1}$ and $T\in\Z$. Then
	\begin{equation*}
		S_{r,w}(T)\geq\frac1w\Bigl(f_{r+1}(T)+\frac{w-1}{2}\,f_r(T)\Bigr).
	\end{equation*}
\end{lem}
\begin{proof}
	Let $\Psi(T)\coloneqq\frac1w\bigl(f_{r+1}(T)+\frac{w-1}{2}\,f_r(T)\bigr)$ be the right-hand
	side of the statement. We prove $S_{r,w}(T)\geq\Psi(T)$ by induction on $T$.

	\emph{Base case $T<0$.} All the arguments of $f_r$ and $f_{r+1}$ in the statement are negative,
	so both sides are $0$ by Definition~\ref{dfn:trunc}.

	\emph{Inductive step.} Let $T\geq0$, and assume the inequality for every integer smaller than
	$T$. We have to prove that $S_{r,w}(T)\geq\Psi(T)$. By the definition of $\Psi$, then by
	Lemma~\ref{lem:increments}~(i) with $u=w$, applied to $f_{r+1}$ and to $f_r$, and finally by
	Lemma~\ref{lem:block-lower} with $t=T$,
	\begin{align*}
		\Psi(T)-\Psi(T-w)
		&=\frac1w\Bigl(f_{r+1}(T)-f_{r+1}(T-w)+\frac{w-1}{2}\bigl(f_r(T)-f_r(T-w)\bigr)\Bigr)\\
		&=\frac1w\Biggl(\,\sum_{u=0}^{w-1}f_r(T-u)+\frac{w-1}{2}\sum_{u=0}^{w-1}f_{r-1}(T-u)\Biggr)
		\leq f_r(T).
	\end{align*}
	Hence, by~\eqref{eq:S-rec}, by the inductive hypothesis for $T-w<T$ and by this inequality,
	\begin{equation*}
		S_{r,w}(T)
		=f_r(T)+S_{r,w}(T-w)
		\geq f_r(T)+\Psi(T-w)
		\geq\Psi(T),
	\end{equation*}
	which is what we had to prove.
\end{proof}

\begin{rem}\label{rem:sum-lower}
	\begin{enumerate}
		\item[(a)] \emph{How much is gained, on a numerical example.} For $r=1$, $w=3$ and $T=7$:
		      \begin{itemize}
			      \item the sum is $S_{1,3}(7)=f_1(7)+f_1(4)+f_1(1)=8+5+2=15$;
			      \item the estimate of Agnarsson, that is Lemma~\ref{lem:sum-lower} without its
			            last term, gives $\frac13\,f_2(7)=\frac{36}{3}=12$;
			      \item Lemma~\ref{lem:sum-lower} gives
			            $\frac13\bigl(f_2(7)+\frac{3-1}{2}\,f_1(7)\bigr)=\frac{36+8}{3}=\frac{44}{3}\approx14.7$.
		      \end{itemize}
		\item[(b)] \emph{The relation with Lemma~3.2 of~\cite{agnarsson}.} That lemma is written
		      with other letters. With $\rho$ in place of the letter $r$ used there, which is not
		      our $r$, its lower bound is
		      \begin{equation*}
			      \sum_{i=0}^{N}\binom{f(i)}{\rho-1}
			      \geq\frac1q\binom{f(0)+1}{\rho},
			      \qquad\text{where }f(i)=P-qi\text{ and }N=\floor*{\frac{P-\rho+1}{q}}.
		      \end{equation*}
		      We take $q=w$, $\rho=r+1$ and $P=T+r$, with $T\in\N$. Then
		      $N=\floor*{\frac{T}{w}}$,
		      $\binom{f(i)}{\rho-1}=\binom{T-iw+r}{r}=f_r(T-iw)$, $\frac1q=\frac1w$ and
		      $\binom{f(0)+1}{\rho}=\binom{T+r+1}{r+1}=f_{r+1}(T)$. So the inequality reads
		      \begin{equation*}
			      S_{r,w}(T)
			      =\sum_{i=0}^{\floor{T/w}}f_r(T-iw)
			      \geq\frac1w\,f_{r+1}(T),
		      \end{equation*}
		      which is Lemma~\ref{lem:sum-lower} without its last term.
		\item[(c)] \emph{The assumption $r\geq1$.} It cannot be dropped. For $r=0$, $T=2$ and $w=3$
		      the sum has the single term $f_0(2)=1$, while the right-hand side would be
		      $\frac13\bigl(f_1(2)+\frac{3-1}{2}f_0(2)\bigr)=\frac{3+1}{3}=\frac43$.
	\end{enumerate}
\end{rem}

\subsection{The theorem}\label{s:lower-thm}
The next theorem adds to the lower bound in~\eqref{eq:agnarsson-ours} binomial coefficients of
lower degree, with the
coefficients of Definition~\ref{dfn:coefficients} for the sequence $b-1$.

\begin{thm}\label{thm:main-lower}
	Let $k\in\N_{\geq2}$ and let $n\in\N$ be a multiple of $d_k$ with $n\geq B_k+k-1$. Then
	\begin{equation*}
		D_k(n)\geq\frac{1}{b_2\dotsm b_k}\sum_{i=0}^{k-2}\bbinom[b-1]{k-2}{i}\binom{n-B_k-i}{k-1-i}.
	\end{equation*}
\end{thm}
\begin{proof}[Proof of Theorem~\ref{thm:main-lower}]
	\emph{The bound as truncated binomial coefficients.} The lemmas used in this proof are stated
	for the truncated binomial coefficients $f_r$, so we first write the statement with them.
	Since $B_k=\Delta_k-(k-1)$ by~\eqref{eq:AB-ours}, the hypothesis $n\geq B_k+k-1$ reads
	$n\geq\Delta_k$. Then $n-\Delta_k\geq0$, so by Definition~\ref{dfn:trunc}, for
	$0\leq i\leq k-2$,
	\begin{equation*}
		f_{k-1-i}(n-\Delta_k)
		=\binom{n-\Delta_k+k-1-i}{k-1-i}
		=\binom{n-B_k-i}{k-1-i}.
	\end{equation*}
	Hence the statement reads
	\begin{equation}\label{eq:claim-lower}
		D_k(n)\geq\frac{1}{b_2\dotsm b_k}\sum_{i=0}^{k-2}\bbinom[b-1]{k-2}{i}f_{k-1-i}(n-\Delta_k)
	\end{equation}
	for every multiple $n\geq\Delta_k$ of $d_k$. If $n$ is a multiple of $d_k$ with $n<\Delta_k$,
	then $n-\Delta_k<0$, so the right-hand side of~\eqref{eq:claim-lower} is $0$ by
	Definition~\ref{dfn:trunc}, and~\eqref{eq:claim-lower} holds because $D_k(n)\geq0$. So it is
	enough to prove~\eqref{eq:claim-lower} for every multiple $n\in\N$ of $d_k$; in this form the
	inductive hypothesis is available for every term of the recursion. We
	prove~\eqref{eq:claim-lower} by induction on $k$.

	\emph{Base case $k=2$.} The sum has the single term $i=0$, with $\bbinom[b-1]{0}{0}=1$,
	so~\eqref{eq:claim-lower} reads $D_2(n)\geq\frac{1}{b_2}f_1(n-\Delta_2)$. For $n<\Delta_2$ it
	has been proved above. Let $n\geq\Delta_2$ be a multiple of $d_2$. By Proposition~3.1
	of~\cite{ballesio-arxiv}, $D_2(n)>\frac{n-\Delta_2}{b_2}$, that is $b_2D_2(n)>n-\Delta_2$. Both
	sides are integers, hence $b_2D_2(n)\geq n-\Delta_2+1$. So, by Definition~\ref{dfn:trunc} and
	since $n-\Delta_2\geq0$,
	\begin{equation*}
		D_2(n)
		\geq\frac{n-\Delta_2+1}{b_2}
		=\frac{1}{b_2}\binom{n-\Delta_2+1}{1}
		=\frac{1}{b_2}\,f_1(n-\Delta_2).
	\end{equation*}

	\emph{Inductive step.} Let $k\geq2$, and assume that, for every multiple $m\in\N$ of $d_k$,
	\begin{equation}\label{eq:lower-hyp}
		D_k(m)\geq\frac{1}{b_2\dotsm b_k}\sum_{i=0}^{k-2}\bbinom[b-1]{k-2}{i}f_{k-1-i}(m-\Delta_k).
	\end{equation}
	We have to prove that, for every multiple $n\in\N$ of $d_{k+1}$,
	\begin{equation}\label{eq:lower-goal}
		D_{k+1}(n)\geq\frac{1}{b_2\dotsm b_{k+1}}\sum_{i=0}^{k-1}\bbinom[b-1]{k-1}{i}f_{k-i}(n-\Delta_{k+1}).
	\end{equation}
	For $n<\Delta_{k+1}$ the right-hand side is $0$, as shown above. So let $n\geq\Delta_{k+1}$ be
	a multiple of $d_{k+1}$, let $s$ be as in Lemma~\ref{lem:recursion}, and set
	$T\coloneqq n-\Delta_{k+1}\in\N$.

	\begin{itemize}
		\item[$\scriptstyle\blacklozenge$] \emph{Step~1: the inductive hypothesis.} By Lemma~\ref{lem:recursion},
		\begin{equation}\label{eq:lower-rec}
			D_{k+1}(n)=\sum_{q=0}^{\floor*{\frac{n-s}{b_{k+1}}}}D_k(n-s-qb_{k+1}),
		\end{equation}
		and every argument $n-s-qb_{k+1}$ with $0\leq q\leq\floor*{\frac{n-s}{b_{k+1}}}$ is a non-negative multiple of $d_k$. So
		we may apply the inductive hypothesis~\eqref{eq:lower-hyp} with $m=n-s-qb_{k+1}$: for
		$0\leq q\leq\floor*{\frac{n-s}{b_{k+1}}}$,
		\begin{equation*}
			D_k(n-s-qb_{k+1})\geq\frac{1}{b_2\dotsm b_k}\sum_{i=0}^{k-2}\bbinom[b-1]{k-2}{i}f_{k-1-i}(n-s-qb_{k+1}-\Delta_k).
		\end{equation*}
		Summing over $0\leq q\leq\floor*{\frac{n-s}{b_{k+1}}}$, by~\eqref{eq:lower-rec} and exchanging the two finite sums,
		\begin{equation*}
			D_{k+1}(n)
			=\sum_{q=0}^{\floor*{\frac{n-s}{b_{k+1}}}}D_k(n-s-qb_{k+1})
			\geq\frac{1}{b_2\dotsm b_k}\sum_{i=0}^{k-2}\bbinom[b-1]{k-2}{i}\sum_{q=0}^{\floor*{\frac{n-s}{b_{k+1}}}}f_{k-1-i}(n-s-qb_{k+1}-\Delta_k).
		\end{equation*}
		The inner sum does not change if it is extended to all $q\geq0$. Indeed, if $q>\floor*{\frac{n-s}{b_{k+1}}}$, then
		$qb_{k+1}>n-s$, that is $n-s-qb_{k+1}<0$, hence $n-s-qb_{k+1}-\Delta_k<0$ because $\Delta_k\geq0$
		(Lemma~\ref{lem:basic}~(iii)), and $f_{k-1-i}$ vanishes at negative integers
		(Definition~\ref{dfn:trunc}). Therefore
		\begin{equation}\label{eq:lower-step1}
			D_{k+1}(n)\geq\frac{1}{b_2\dotsm b_k}\sum_{i=0}^{k-2}\bbinom[b-1]{k-2}{i}\sum_{q\geq0}f_{k-1-i}(n-s-qb_{k+1}-\Delta_k).
		\end{equation}

		\item[$\scriptstyle\blacklozenge$] \emph{Step~2: comparison with a progression.} By Lemma~\ref{lem:basic}~(iii),
		$\Delta_{k+1}=\Delta_k+b_{k+1}-a_{k+1}$, so $T=n-(b_{k+1}-a_{k+1})-\Delta_k$; and
		$s\leq b_{k+1}-a_{k+1}$ by Lemma~\ref{lem:recursion}. Hence, for every $q\geq0$,
		\begin{equation*}
			n-s-qb_{k+1}-\Delta_k
			\geq n-(b_{k+1}-a_{k+1})-qb_{k+1}-\Delta_k
			=T-qb_{k+1}.
		\end{equation*}
		Since $f_{k-1-i}$ is non-decreasing (Lemma~\ref{lem:trunc}~(ii)), summing over $q\geq0$ and
		using~\eqref{eq:S-def},
		\begin{equation*}
			\sum_{q\geq0}f_{k-1-i}(n-s-qb_{k+1}-\Delta_k)
			\geq\sum_{q\geq0}f_{k-1-i}(T-qb_{k+1})
			=S_{k-1-i,\,b_{k+1}}(T).
		\end{equation*}
		As the coefficients $\bbinom[b-1]{k-2}{i}$ are non-negative (Lemma~\ref{lem:sym-poly}~(ii)),
		from~\eqref{eq:lower-step1} we get
		\begin{equation}\label{eq:lower-step2}
			D_{k+1}(n)\geq\frac{1}{b_2\dotsm b_k}\sum_{i=0}^{k-2}\bbinom[b-1]{k-2}{i}S_{k-1-i,\,b_{k+1}}(T).
		\end{equation}

		\item[$\scriptstyle\blacklozenge$] \emph{Step~3: the sums along the progression.} For $0\leq i\leq k-2$ we have $k-1-i\geq1$, and
		Lemma~\ref{lem:sum-lower} with $r=k-1-i$ and $w=b_{k+1}$ gives
		\begin{equation*}
			S_{k-1-i,\,b_{k+1}}(T)
			\geq\frac{1}{b_{k+1}}\Bigl(f_{k-i}(T)+\frac{b_{k+1}-1}{2}\,f_{k-1-i}(T)\Bigr).
		\end{equation*}
		Since the coefficients $\bbinom[b-1]{k-2}{i}$ are non-negative, we may insert this estimate
		in~\eqref{eq:lower-step2}; then, using $b_2\dotsm b_k\cdot b_{k+1}=b_2\dotsm b_{k+1}$,
		\begin{align}
			D_{k+1}(n)
			&\geq\frac{1}{b_2\dotsm b_k}\sum_{i=0}^{k-2}\bbinom[b-1]{k-2}{i}\frac{1}{b_{k+1}}\Bigl(f_{k-i}(T)+\frac{b_{k+1}-1}{2}\,f_{k-1-i}(T)\Bigr)\notag\\
			&=\frac{1}{b_2\dotsm b_{k+1}}\Biggl(\,\sum_{i=0}^{k-2}\bbinom[b-1]{k-2}{i}f_{k-i}(T)
			+\frac{b_{k+1}-1}{2}\sum_{i=0}^{k-2}\bbinom[b-1]{k-2}{i}f_{k-1-i}(T)\Biggr).
			\label{eq:lower-step3}
		\end{align}

		\item[$\scriptstyle\blacklozenge$] \emph{Step~4: collecting the coefficients.} The expression in brackets
		in~\eqref{eq:lower-step3} is the left-hand side of Lemma~\ref{lem:pascal-sums} with $c=b-1$,
		$m=k-1$ (so that $c_{m+2}=b_{k+1}-1$) and $y_i\coloneqq f_{k-i}(T)$ for $0\leq i\leq k-1$;
		indeed $y_{i+1}=f_{k-1-i}(T)$. Hence, by~\eqref{eq:lower-step3}, then by
		Lemma~\ref{lem:pascal-sums}, and finally since $T=n-\Delta_{k+1}$,
		\begin{align*}
			D_{k+1}(n)
			&\geq\frac{1}{b_2\dotsm b_{k+1}}\Biggl(\,\sum_{i=0}^{k-2}\bbinom[b-1]{k-2}{i}f_{k-i}(T)
			+\frac{b_{k+1}-1}{2}\sum_{i=0}^{k-2}\bbinom[b-1]{k-2}{i}f_{k-1-i}(T)\Biggr)\\
			&=\frac{1}{b_2\dotsm b_{k+1}}\sum_{i=0}^{k-1}\bbinom[b-1]{k-1}{i}f_{k-i}(T)
			=\frac{1}{b_2\dotsm b_{k+1}}\sum_{i=0}^{k-1}\bbinom[b-1]{k-1}{i}f_{k-i}(n-\Delta_{k+1}),
		\end{align*}
		which is~\eqref{eq:lower-goal}.
	\end{itemize}
\end{proof}

\subsection{Comparison with the lower bound of Agnarsson}\label{s:lower-agnarsson}
Both the lower bound of Theorem~\ref{thm:main-lower} and that of
Agnarsson in~\eqref{eq:agnarsson-ours} have the factor $\frac{1}{b_2\dotsm b_k}$, so in the
comparison we leave it out.

\begin{prop}[Agnarsson's lower bound and Theorem~\ref{thm:main-lower}]\label{prop:agnarsson-lower}
	Let $k\in\N_{\geq2}$. For every integer $n\geq B_k+k-1$, the lower bound of
	Theorem~\ref{thm:main-lower} is greater than or equal to the lower bound of Agnarsson. It is
	greater if $k\geq3$ and $b_i\geq2$ for some $3\leq i\leq k$.
\end{prop}
\begin{proof}
	Let $n\geq B_k+k-1$. We have to prove that
	\begin{equation*}
		\sum_{i=0}^{k-2}\bbinom[b-1]{k-2}{i}\binom{n-B_k-i}{k-1-i}\geq\binom{n-B_k}{k-1},
	\end{equation*}
	with strict inequality if $k\geq3$ and $b_i\geq2$ for some $3\leq i\leq k$. The term $i=0$ of
	the sum is $\binom{n-B_k}{k-1}$, since $\bbinom[b-1]{k-2}{0}=1$
	(Lemma~\ref{lem:sym-poly}~(i)). For $1\leq i\leq k-2$ the coefficient $\bbinom[b-1]{k-2}{i}$
	is non-negative (Lemma~\ref{lem:sym-poly}~(ii)), and $n-B_k-i\geq k-1-i\geq0$, so the
	binomial coefficient $\binom{n-B_k-i}{k-1-i}$ is a positive integer. Hence
	\begin{equation*}
		\sum_{i=0}^{k-2}\bbinom[b-1]{k-2}{i}\binom{n-B_k-i}{k-1-i}
		=\binom{n-B_k}{k-1}+\sum_{i=1}^{k-2}\bbinom[b-1]{k-2}{i}\binom{n-B_k-i}{k-1-i}
		\geq\binom{n-B_k}{k-1}.
	\end{equation*}
	If $k\geq3$ and $b_i\geq2$ for some $3\leq i\leq k$, then
	$\bbinom[b-1]{k-2}{1}=\sum_{l=3}^{k}\frac{b_l-1}{2}>0$, so the term $i=1$ is positive and the
	inequality is strict.
\end{proof}

\begin{rem}[The assumption $D_k(n)>0$]\label{rem:positivity}
	Agnarsson's proof of the lower bound starts from a solution of $a_1x_1+\dotsb+a_kx_k=n$: from
	it he derives the recursion for $D_k(n)$ that drives the induction (formula~(4)
	of~\cite{agnarsson}, based on his Observation~3.1). This is why he assumes $D_k(n)>0$. Our proofs
	use instead the recursion of Lemma~\ref{lem:recursion}, which needs only that $n$ is a
	multiple of $d_k$: the integer $s$ that starts the arithmetic progression is obtained from
	B\'ezout's identity, not from a solution. Hence Theorem~\ref{thm:main-lower} holds without the
	assumption $D_k(n)>0$, and it gives $D_k(n)\geq\frac{1}{b_2\dotsm b_k}>0$ for every multiple
	$n\geq B_k+k-1$ of $d_k$, since its term $i=0$ is $\frac{1}{b_2\dotsm b_k}\binom{n-B_k}{k-1}$,
	with $n-B_k\geq k-1$.
\end{rem}

\section{The upper bound}\label{s:upper}
Agnarsson's upper bound in~\eqref{eq:agnarsson-ours} is a single binomial coefficient.
Theorem~\ref{thm:main-upper} improves it: it replaces $A_k$ by the number $M_k$
of~\eqref{eq:M-def}, which is not larger, and it adds terms of lower degree when some of the
numbers $b_3,\dotsc,b_k$ are even. Its proof needs an upper estimate for the
sums~\eqref{eq:S-def} and some averages of truncated binomial coefficients, which we introduce
first.

\subsection{An upper estimate for the sums along a progression}\label{s:upper-sums}
For the upper estimate we use blocks centred at $t$. For $w\in\N_{\geq1}$ let
\begin{equation}\label{eq:UV}
	U\coloneqq\floor*{\frac{w}{2}},\qquad V\coloneqq\floor*{\frac{w-1}{2}},
\end{equation}
so that $U+V+1=w$, and $V=U$ if $w$ is odd, $V=U-1$ if $w$ is even. The block
$t-V,\dotsc,t+U$ and its mirror image $t-U,\dotsc,t+V$ have $w$ elements each. If $w$ is odd they
coincide; if $w$ is even, the first has one more element above $t$ than below, and the second one
more below than above.

\begin{lem}\label{lem:block-upper}
	Let $r,w\in\N_{\geq1}$ and $t\in\Z$, and let $U$ and $V$ be as in~\eqref{eq:UV}. Then
	\begin{equation*}
		w\,f_r(t)\leq\frac12\sum_{u=-V}^{U}f_r(t+u)+\frac12\sum_{v=-U}^{V}f_r(t+v).
	\end{equation*}
\end{lem}
For instance, for $w=3$ the right-hand side is $f_r(t-1)+f_r(t)+f_r(t+1)$, and for $w=4$ it is
$\frac12f_r(t-2)+f_r(t-1)+f_r(t)+f_r(t+1)+\frac12f_r(t+2)$.
\begin{proof}
	Putting $v=-u$ in the second sum, it becomes $\sum_{u=-V}^{U}f_r(t-u)$. By
	Lemma~\ref{lem:increments}~(ii) with $j=\abs{u}$, we have $f_r(t+u)+f_r(t-u)\geq2f_r(t)$ for
	every $u$; and the sum has $U+V+1=w$ terms. Hence
	\begin{equation*}
		\frac12\sum_{u=-V}^{U}f_r(t+u)+\frac12\sum_{v=-U}^{V}f_r(t+v)
		=\frac12\sum_{u=-V}^{U}\bigl(f_r(t+u)+f_r(t-u)\bigr)
		\geq\frac12\sum_{u=-V}^{U}2f_r(t)
		=w\,f_r(t).\qedhere
	\end{equation*}
\end{proof}

\begin{lem}\label{lem:sum-upper}
	Let $r,w\in\N_{\geq1}$ and $T\in\Z$, and let $U$ and $V$ be as in~\eqref{eq:UV}. Then
	\begin{equation*}
		S_{r,w}(T)\leq\frac{1}{2w}\Bigl(f_{r+1}(T+U)+f_{r+1}(T+V)\Bigr).
	\end{equation*}
	In particular, if $w$ is odd then $S_{r,w}(T)\leq\frac1w f_{r+1}\bigl(T+\frac{w-1}{2}\bigr)$.
\end{lem}
\begin{proof}
	Let $\Phi(T)\coloneqq\frac{1}{2w}\bigl(f_{r+1}(T+U)+f_{r+1}(T+V)\bigr)$ be the right-hand
	side of the statement. We prove $S_{r,w}(T)\leq\Phi(T)$ by induction on $T$.

	\emph{Base case $T<0$.} All the arguments of $f_r$ in $S_{r,w}(T)$ are negative, so
	$S_{r,w}(T)=0$ by Definition~\ref{dfn:trunc}, and $\Phi(T)\geq0$ because $f_{r+1}$ is
	non-negative.

	\emph{Inductive step.} Let $T\geq0$, and assume the inequality for every integer smaller than
	$T$. We have to prove that $S_{r,w}(T)\leq\Phi(T)$. By the definition of $\Phi$, then by
	Lemma~\ref{lem:increments}~(i) for the increments of $f_{r+1}$ over $w$ steps at the points
	$T+U$ and $T+V$, then putting $u=U-j$ and $v=V-i$, where $U-(w-1)=-V$ and $V-(w-1)=-U$, and
	finally by
	Lemma~\ref{lem:block-upper} with $t=T$,
	\begin{align*}
		\Phi(T)-\Phi(T-w)
		&=\frac{1}{2w}\Bigl(f_{r+1}(T+U)-f_{r+1}(T+U-w)+f_{r+1}(T+V)-f_{r+1}(T+V-w)\Bigr)\\
		&=\frac{1}{2w}\Biggl(\,\sum_{j=0}^{w-1}f_r(T+U-j)+\sum_{i=0}^{w-1}f_r(T+V-i)\Biggr)\\
		&=\frac{1}{2w}\Biggl(\,\sum_{u=-V}^{U}f_r(T+u)+\sum_{v=-U}^{V}f_r(T+v)\Biggr)
		\geq f_r(T).
	\end{align*}
	Hence, by~\eqref{eq:S-rec}, by the inductive hypothesis for $T-w<T$ and by this inequality,
	\begin{equation*}
		S_{r,w}(T)
		=f_r(T)+S_{r,w}(T-w)
		\leq f_r(T)+\Phi(T-w)
		\leq\Phi(T),
	\end{equation*}
	which is what we had to prove.
\end{proof}

\begin{rem}\label{rem:sum-upper}
	\begin{enumerate}
		\item[(a)] \emph{How much is gained, on a numerical example.} For $r=1$, $w=3$ and $T=7$:
		      \begin{itemize}
			      \item the sum is $S_{1,3}(7)=f_1(7)+f_1(4)+f_1(1)=8+5+2=15$;
			      \item the estimate of Agnarsson, recalled in~(b), gives
			            $\frac13\,f_2(9)=\frac{55}{3}\approx18.3$;
			      \item Lemma~\ref{lem:sum-upper}, where $U=V=1$, gives
			            $\frac16\bigl(f_2(8)+f_2(8)\bigr)=\frac{45+45}{6}=15$.
		      \end{itemize}
		      For an even difference, $w=4$, with the same $r$ and $T$: the sum is
		      $S_{1,4}(7)=f_1(7)+f_1(3)=8+4=12$; the estimate of Agnarsson gives
		      $\frac14\,f_2(10)=\frac{66}{4}=16.5$; and Lemma~\ref{lem:sum-upper}, where $U=2$ and
		      $V=1$, gives $\frac18\bigl(f_2(9)+f_2(8)\bigr)=\frac{55+45}{8}=12.5$.
		\item[(b)] \emph{The relation with Lemma~3.2 of~\cite{agnarsson}.} With the letters of
		      Remark~\ref{rem:sum-lower}~(b), the upper bound of that lemma is
		      \begin{equation*}
			      \sum_{i=0}^{N}\binom{f(i)}{\rho-1}
			      \leq\frac1q\binom{f(-1)}{\rho},
			      \qquad\text{where }f(i)=P-qi\text{ and }N=\floor*{\frac{P-\rho+1}{q}}.
		      \end{equation*}
		      We take again $q=w$, $\rho=r+1$ and $P=T+r$, with $T\in\N$. Then
		      $N=\floor*{\frac{T}{w}}$,
		      $\binom{f(i)}{\rho-1}=f_r(T-iw)$, $\frac1q=\frac1w$ and, since $f(-1)=P+q=T+r+w$,
		      $\binom{f(-1)}{\rho}=\binom{(T+w-1)+r+1}{r+1}=f_{r+1}(T+w-1)$. So the inequality reads
		      \begin{equation*}
			      S_{r,w}(T)
			      =\sum_{i=0}^{\floor{T/w}}f_r(T-iw)
			      \leq\frac1w\,f_{r+1}(T+w-1)
			      =\frac{1}{2w}\Bigl(f_{r+1}(T+w-1)+f_{r+1}(T+w-1)\Bigr).
		      \end{equation*}
		      In the last form it can be compared term by term with the bound
		      $\frac{1}{2w}\bigl(f_{r+1}(T+U)+f_{r+1}(T+V)\bigr)$ of Lemma~\ref{lem:sum-upper},
		      which is not larger: $V\leq U\leq w-1$ and $f_{r+1}$ is non-decreasing
		      (Lemma~\ref{lem:trunc}~(ii)), so both $f_{r+1}(T+U)$ and $f_{r+1}(T+V)$ are at most
		      $f_{r+1}(T+w-1)$.
		\item[(c)] \emph{The assumption $r\geq1$.} It cannot be dropped. For $r=0$, $T=0$ and $w=3$
		      the sum has the single term $f_0(0)=1$, while the right-hand side, where $U=V=1$,
		      would be $\frac16\bigl(f_1(1)+f_1(1)\bigr)=\frac{2+2}{6}=\frac23$.
	\end{enumerate}
\end{rem}

\subsection{Parities and averages}\label{s:averages}
The upper bound depends on the parities of $b_3,\dotsc,b_k$. For $i\geq2$
and $k\geq2$ let
\begin{equation}\label{eq:eps-E-def}
	\varepsilon_i\coloneqq
	\begin{cases}
		1 & \text{if $b_i$ is even},\\
		0 & \text{if $b_i$ is odd},
	\end{cases}
	\qquad
	E_k\coloneqq\sum_{i=3}^{k}\varepsilon_i,
\end{equation}
so that $E_k$ is the number of even numbers among $b_3,\dotsc,b_k$. By these definitions and
the definition~\eqref{eq:M-def} of $M_k$, we have $M_2=b_2$, $E_2=0$ and, for $k\geq2$,
\begin{equation}\label{eq:M-rec}
	M_{k+1}=M_k+\floor*{\frac{b_{k+1}}{2}}+1,\qquad E_{k+1}=E_k+\varepsilon_{k+1},
	\qquad \floor*{\frac{b_{k+1}-1}{2}}=\floor*{\frac{b_{k+1}}{2}}-\varepsilon_{k+1}.
\end{equation}

The coefficients of the terms of lower degree of the upper bound are given by
Definition~\ref{dfn:coefficients} with the sequence $c_l=-\varepsilon_l$, which we denote by
$-\varepsilon$:
\begin{equation}\label{eq:coeff-eps}
	\bbinom[-\varepsilon]{m}{i}
	=e_i\Bigl(-\frac{\varepsilon_3}{2},\dotsc,-\frac{\varepsilon_{m+2}}{2}\Bigr)
	=\Bigl(-\frac12\Bigr)^{i}\binom{E_{m+2}}{i}\qquad(0\leq i\leq m),
\end{equation}
since exactly $E_{m+2}$ of the numbers $-\frac{\varepsilon_3}{2},\dotsc,-\frac{\varepsilon_{m+2}}{2}$
are equal to $-\frac12$ and the others are $0$.

These coefficients have alternating sign, and an upper estimate cannot be applied term by term to a
sum with such coefficients, as is done for the lower bound in the proof of Theorem~\ref{thm:main-lower}. For this
reason we use averages of binomial coefficients of degree $k-1$ with non-negative weights. The next
lemma introduces these averages and relates them to the sums with the coefficients~\eqref{eq:coeff-eps}.

\begin{lem}\label{lem:average}
	For $k\in\N_{\geq2}$, $r\in\N$ and $t\in\Z$ let
	\begin{equation}\label{eq:average-def}
		\alpha_{k,r}(t)\coloneqq\frac{1}{2^{E_k}}\sum_{j=0}^{E_k}\binom{E_k}{j}f_r(t-j).
	\end{equation}
	Then:
	\begin{enumerate}
		\item[\textup{(i)}] $\alpha_{2,r}=f_r$, and
		      $\alpha_{k+1,r}(t)=\frac12\bigl(\alpha_{k,r}(t)+\alpha_{k,r}(t-\varepsilon_{k+1})\bigr)$;
		\item[\textup{(ii)}] $\alpha_{k,r}$ is non-negative and non-decreasing, and
		      $\alpha_{k,r}(t)\leq f_r(t)$;
		\item[\textup{(iii)}] if $r\geq k-1$, then
		      $\displaystyle\alpha_{k,r}(t)=\sum_{i=0}^{k-2}\bbinom[-\varepsilon]{k-2}{i}f_{r-i}(t)$.
	\end{enumerate}
\end{lem}
\begin{proof}
	\begin{enumerate}
		\item[\textup{(i)}] $E_2=0$, so $\alpha_{2,r}=f_r$. If $\varepsilon_{k+1}=0$, then
		      $E_{k+1}=E_k$, so $\alpha_{k+1,r}=\alpha_{k,r}$, and the right-hand side is
		      $\frac12\bigl(\alpha_{k,r}(t)+\alpha_{k,r}(t)\bigr)=\alpha_{k,r}(t)$. If
		      $\varepsilon_{k+1}=1$, then $E_{k+1}=E_k+1$ and
		      \begin{equation*}
		      	\frac12\bigl(\alpha_{k,r}(t)+\alpha_{k,r}(t-1)\bigr)
		      	=\frac{1}{2^{E_k+1}}\sum_{j=0}^{E_k}\binom{E_k}{j}\bigl(f_r(t-j)+f_r(t-j-1)\bigr)
		      	=\frac{1}{2^{E_k+1}}\sum_{j=0}^{E_k+1}\binom{E_k+1}{j}f_r(t-j)
		      	=\alpha_{k+1,r}(t),
		      \end{equation*}
		      because the coefficient of $f_r(t-j)$ is $\binom{E_k}{j}+\binom{E_k}{j-1}=\binom{E_k+1}{j}$
		      by Pascal's rule (with $\binom{E_k}{-1}=\binom{E_k}{E_k+1}=0$).
		\item[\textup{(ii)}] The weights $\frac{1}{2^{E_k}}\binom{E_k}{j}$ are non-negative and
		      their sum is $1$, by the binomial theorem. The functions $f_r(t-j)$ are non-negative and
		      non-decreasing in $t$ (Lemma~\ref{lem:trunc}~(ii)), hence so is $\alpha_{k,r}$; and
		      $f_r(t-j)\leq f_r(t)$ for $j\geq0$, hence
		      $\alpha_{k,r}(t)\leq\frac{1}{2^{E_k}}\sum_{j=0}^{E_k}\binom{E_k}{j}f_r(t)=f_r(t)$.
		\item[\textup{(iii)}] We prove the statement by induction on $k$.

		      \emph{Base case $k=2$.} For $k=2$ the sum has the single term $i=0$, with
		      $\bbinom[-\varepsilon]{0}{0}=1$, so the right-hand side is $f_r(t)$; and
		      $\alpha_{2,r}(t)=f_r(t)$ by~(i).

		      \emph{Inductive step.} Assume the statement for some $k\geq2$ and every $r\geq k-1$.
		      We have to prove it for $k+1$ and every $r\geq k$; let $r\geq k$. By
		      Lemma~\ref{lem:trunc}~(i), applied to each term
		      of~\eqref{eq:average-def},
		      \begin{equation*}
		      	\alpha_{k,r}(t)-\alpha_{k,r}(t-1)
		      	=\frac{1}{2^{E_k}}\sum_{j=0}^{E_k}\binom{E_k}{j}\bigl(f_r(t-j)-f_r(t-j-1)\bigr)
		      	=\frac{1}{2^{E_k}}\sum_{j=0}^{E_k}\binom{E_k}{j}f_{r-1}(t-j)
		      	=\alpha_{k,r-1}(t).
		      \end{equation*}
		      Hence, by~(i), if $\varepsilon_{k+1}=1$,
		      \begin{equation*}
		      	\alpha_{k+1,r}(t)
		      	=\frac12\bigl(\alpha_{k,r}(t)+\alpha_{k,r}(t-1)\bigr)
		      	=\alpha_{k,r}(t)-\frac12\bigl(\alpha_{k,r}(t)-\alpha_{k,r}(t-1)\bigr)
		      	=\alpha_{k,r}(t)-\frac12\,\alpha_{k,r-1}(t),
		      \end{equation*}
		      while $\alpha_{k+1,r}(t)=\alpha_{k,r}(t)$ if $\varepsilon_{k+1}=0$. In both cases
		      $\alpha_{k+1,r}(t)=\alpha_{k,r}(t)-\frac{\varepsilon_{k+1}}{2}\,\alpha_{k,r-1}(t)$, and
		      the inductive hypothesis for $r$ and for $r-1\geq k-1$ gives
		      \begin{equation*}
		      	\alpha_{k+1,r}(t)
		      	=\sum_{i=0}^{k-2}\bbinom[-\varepsilon]{k-2}{i}f_{r-i}(t)
		      	-\frac{\varepsilon_{k+1}}{2}\sum_{i=0}^{k-2}\bbinom[-\varepsilon]{k-2}{i}f_{r-1-i}(t).
		      \end{equation*}
		      The right-hand side is the left-hand side of Lemma~\ref{lem:pascal-sums} with
		      $c=-\varepsilon$, $m=k-1$ (so that $c_{m+2}=-\varepsilon_{k+1}$) and
		      $y_i\coloneqq f_{r-i}(t)$ for $0\leq i\leq k-1$; indeed $y_{i+1}=f_{r-1-i}(t)$. Hence
		      \begin{equation*}
		      	\alpha_{k+1,r}(t)=\sum_{i=0}^{k-1}\bbinom[-\varepsilon]{k-1}{i}f_{r-i}(t),
		      \end{equation*}
		      which is the statement for $k+1$.\qedhere
	\end{enumerate}
\end{proof}

\subsection{The theorem}\label{s:upper-thm}
The next theorem replaces $A_k$ in~\eqref{eq:agnarsson-ours} by the number $M_k$ defined
in~\eqref{eq:M-def}, and adds terms of lower degree with the coefficients~\eqref{eq:coeff-eps}.

\begin{thm}\label{thm:main-upper}
	Let $k\in\N_{\geq2}$. Then for every $n\in\N$
	\begin{align}
		D_k(n)
		&\leq\frac{1}{b_2\dotsm b_k}\sum_{i=0}^{k-2}\bbinom[-\varepsilon]{k-2}{i}\binom{n+M_k-i}{k-1-i}
		\label{eq:upper-sum}\\
		&=\frac{1}{b_2\dotsm b_k}\cdot\frac{1}{2^{E_k}}\sum_{j=0}^{E_k}\binom{E_k}{j}\binom{n+M_k-j}{k-1}
		\label{eq:upper-average}\\
		&\leq\frac{1}{b_2\dotsm b_k}\binom{n+M_k}{k-1}.
		\label{eq:upper-single}
	\end{align}
	If $b_3,\dotsc,b_k$ are all odd, the last inequality is an equality; otherwise it is strict.
\end{thm}
\begin{proof}[Proof of Theorem~\ref{thm:main-upper}]
	\emph{The bounds as truncated binomial coefficients.} The lemmas used in this proof are stated
	for the truncated binomial coefficients $f_r$, so we first write the statement with them.
	By Definition~\ref{dfn:trunc}, $f_r(t)=\binom{t+r}{r}$ only for $t\geq0$; so we have to check
	that the arguments are non-negative. Throughout the proof we set
	\begin{equation*}
		T\coloneqq n+M_k-(k-1).
	\end{equation*}
	The binomial coefficients of the statement are
	$\binom{n+M_k-i}{k-1-i}=\binom{T+k-1-i}{k-1-i}$ in~\eqref{eq:upper-sum}, for $0\leq i\leq k-2$,
	and $\binom{n+M_k-j}{k-1}=\binom{(T-j)+k-1}{k-1}$ in~\eqref{eq:upper-average}, for
	$0\leq j\leq E_k$; recall from~\eqref{eq:eps-E-def} that $E_k$ is the number of
	even numbers among $b_3,\dotsc,b_k$. So the arguments to be checked are $T$ and $T-j$ for
	$0\leq j\leq E_k$, and the smallest of them is $T-E_k$.

	To bound $T-E_k$ from below we compare $M_k$ with $E_k$ term by term. For every $i\geq3$ we have
	$\floor*{\frac{b_i}{2}}\geq\varepsilon_i$: if $b_i$ is even, then $b_i\geq2$ and
	$\floor*{\frac{b_i}{2}}\geq1=\varepsilon_i$; if $b_i$ is odd, then
	$\floor*{\frac{b_i}{2}}\geq0=\varepsilon_i$. Moreover $b_2\geq1$. Hence, by~\eqref{eq:M-def}
	and~\eqref{eq:eps-E-def},
	\begin{equation}\label{eq:M-lower}
		M_k
		=b_2+\sum_{i=3}^{k}\Bigl(\floor*{\frac{b_i}{2}}+1\Bigr)
		\geq1+\sum_{i=3}^{k}(\varepsilon_i+1)
		=k-1+E_k.
	\end{equation}
	Therefore, for $0\leq j\leq E_k$, by~\eqref{eq:M-lower} and since $n\geq0$,
	\begin{equation*}
		T-j
		\geq T-E_k
		=n+M_k-(k-1)-E_k
		\geq n
		\geq0;
	\end{equation*}
	in particular $T\geq0$ (take $j=0$). Unlike in Theorem~\ref{thm:main-lower}, where the
	corresponding condition $n\geq\Delta_k$ is a hypothesis, here non-negativity holds for every
	$n\in\N$. Hence, by Definition~\ref{dfn:trunc}, for $0\leq i\leq k-2$ and $0\leq j\leq E_k$,
	\begin{equation*}
		f_{k-1-i}(T)
		=\binom{T+k-1-i}{k-1-i}
		=\binom{n+M_k-i}{k-1-i},
		\qquad
		f_{k-1}(T-j)
		=\binom{T-j+k-1}{k-1}
		=\binom{n+M_k-j}{k-1}.
	\end{equation*}

	\emph{The equality of~\eqref{eq:upper-sum} and~\eqref{eq:upper-average}.} By the second identity
	and~\eqref{eq:average-def}, and then by Lemma~\ref{lem:average}~(iii) with $r=k-1$ and by the
	first identity,
	\begin{equation*}
		\frac{1}{2^{E_k}}\sum_{j=0}^{E_k}\binom{E_k}{j}\binom{n+M_k-j}{k-1}
		=\alpha_{k,k-1}(T)
		=\sum_{i=0}^{k-2}\bbinom[-\varepsilon]{k-2}{i}\binom{n+M_k-i}{k-1-i}.
	\end{equation*}

	\emph{The inequality~\eqref{eq:upper-single}.} By Lemma~\ref{lem:average}~(ii) and the second
	identity with $j=0$,
	\begin{equation*}
		\alpha_{k,k-1}(T)
		\leq f_{k-1}(T)
		=\binom{n+M_k}{k-1}.
	\end{equation*}

	\emph{Case $b_3,\dotsc,b_k$ all odd.} Then $E_k=0$, so $\alpha_{k,k-1}(T)=f_{k-1}(T)$
	by~\eqref{eq:average-def}, and~\eqref{eq:upper-single} is an equality.

	\emph{Case some $b_i$ even, $3\leq i\leq k$.} Then $E_k\geq1$, hence $T\geq E_k\geq1$. Since the
	weights $\frac{1}{2^{E_k}}\binom{E_k}{j}$ have sum $1$, and $f_{k-1}(T)-f_{k-1}(T-j)\geq0$ for
	every $j\geq0$ (Lemma~\ref{lem:trunc}~(ii)), keeping only the term $j=1$ we get
	\begin{align*}
		\alpha_{k,k-1}(T)
		&=f_{k-1}(T)-\frac{1}{2^{E_k}}\sum_{j=0}^{E_k}\binom{E_k}{j}\bigl(f_{k-1}(T)-f_{k-1}(T-j)\bigr)\\
		&\leq f_{k-1}(T)-\frac{E_k}{2^{E_k}}\bigl(f_{k-1}(T)-f_{k-1}(T-1)\bigr)
		=f_{k-1}(T)-\frac{E_k}{2^{E_k}}\,f_{k-2}(T),
	\end{align*}
	by Lemma~\ref{lem:trunc}~(i). As $f_{k-2}(T)=\binom{T+k-2}{k-2}\geq1$, it follows that
	$\alpha_{k,k-1}(T)<f_{k-1}(T)$, and~\eqref{eq:upper-single} is strict.

	\emph{The bound for $D_k(n)$.} By the equality proved above, the first inequality of the statement
	reads
	\begin{equation}\label{eq:claim-upper}
		D_k(n)\leq\frac{1}{b_2\dotsm b_k}\,\alpha_{k,k-1}\bigl(n+M_k-(k-1)\bigr).
	\end{equation}
	If $n$ is not a multiple of $d_k$, then $D_k(n)=0$ by Lemma~\ref{lem:basic}~(iv),
	and~\eqref{eq:claim-upper} holds because its right-hand side is non-negative
	(Lemma~\ref{lem:average}~(ii)). So it is enough to prove~\eqref{eq:claim-upper} for every
	multiple $n\in\N$ of $d_k$. We prove~\eqref{eq:claim-upper} by induction on $k$.

	\emph{Base case $k=2$.} Here $\alpha_{2,1}=f_1$ by Lemma~\ref{lem:average}~(i) and $M_2=b_2$,
	so~\eqref{eq:claim-upper} reads $D_2(n)\leq\frac{1}{b_2}f_1(n+b_2-1)$. Let $n\in\N$ be a multiple
	of $d_2$. By Proposition~3.1 of~\cite{ballesio-arxiv} (where $b_2$ is written $\nabla_2$),
	$D_2(n)\leq\frac{n+b_2}{b_2}$. So, by Definition~\ref{dfn:trunc} and since $n+b_2-1\geq0$,
	\begin{equation*}
		D_2(n)
		\leq\frac{n+b_2}{b_2}
		=\frac{1}{b_2}\binom{n+b_2}{1}
		=\frac{1}{b_2}\,f_1(n+b_2-1).
	\end{equation*}

	\emph{Inductive step.} Let $k\geq2$, and assume that, for every multiple $m\in\N$ of $d_k$,
	\begin{equation}\label{eq:upper-hyp}
		D_k(m)\leq\frac{1}{b_2\dotsm b_k}\,\alpha_{k,k-1}\bigl(m+M_k-(k-1)\bigr).
	\end{equation}
	We have to prove that, for every multiple $n\in\N$ of $d_{k+1}$,
	\begin{equation}\label{eq:upper-goal}
		D_{k+1}(n)\leq\frac{1}{b_2\dotsm b_{k+1}}\,\alpha_{k+1,k}\bigl(n+M_{k+1}-k\bigr).
	\end{equation}
	Let $n\in\N$ be a multiple of $d_{k+1}$, let $s$ be as in Lemma~\ref{lem:recursion}, and
	let $T=n+M_k-(k-1)$, as above.

	\begin{itemize}
		\item[$\scriptstyle\blacklozenge$] \emph{Step~1: the inductive hypothesis.} By Lemma~\ref{lem:recursion},
		\begin{equation}\label{eq:upper-rec}
			D_{k+1}(n)=\sum_{q=0}^{\floor*{\frac{n-s}{b_{k+1}}}}D_k(n-s-qb_{k+1}),
		\end{equation}
		and every argument $n-s-qb_{k+1}$ with $0\leq q\leq\floor*{\frac{n-s}{b_{k+1}}}$ is a non-negative
		multiple of $d_k$. So we may apply the inductive hypothesis~\eqref{eq:upper-hyp} with
		$m=n-s-qb_{k+1}$: for $0\leq q\leq\floor*{\frac{n-s}{b_{k+1}}}$,
		\begin{equation*}
			D_k(n-s-qb_{k+1})
			\leq\frac{1}{b_2\dotsm b_k}\,\alpha_{k,k-1}(n-s-qb_{k+1}+M_k-(k-1))
			=\frac{1}{b_2\dotsm b_k}\,\alpha_{k,k-1}(T-s-qb_{k+1}).
		\end{equation*}
		Summing over $0\leq q\leq\floor*{\frac{n-s}{b_{k+1}}}$, by~\eqref{eq:upper-rec},
		\begin{equation*}
			D_{k+1}(n)
			=\sum_{q=0}^{\floor*{\frac{n-s}{b_{k+1}}}}D_k(n-s-qb_{k+1})
			\leq\frac{1}{b_2\dotsm b_k}\sum_{q=0}^{\floor*{\frac{n-s}{b_{k+1}}}}\alpha_{k,k-1}(T-s-qb_{k+1}).
		\end{equation*}
		Unlike in the proof of Theorem~\ref{thm:main-lower}, the terms with
		$q>\floor*{\frac{n-s}{b_{k+1}}}$ need not vanish, since $T-s-qb_{k+1}$ may still be non-negative;
		but they are non-negative (Lemma~\ref{lem:average}~(ii)), so adding them can only increase the
		right-hand side. Therefore
		\begin{equation}\label{eq:upper-step1}
			D_{k+1}(n)\leq\frac{1}{b_2\dotsm b_k}\sum_{q\geq0}\alpha_{k,k-1}(T-s-qb_{k+1}),
		\end{equation}
		where only finitely many terms are different from $0$.

		\item[$\scriptstyle\blacklozenge$] \emph{Step~2: comparison with a progression.} Since $s\geq0$, for every $q\geq0$ we have
		$T-s-qb_{k+1}\leq T-qb_{k+1}$, and $\alpha_{k,k-1}$ is non-decreasing
		(Lemma~\ref{lem:average}~(ii)). Hence, by~\eqref{eq:average-def}, exchanging the two sums
		(only finitely many terms are different from $0$) and by~\eqref{eq:S-def},
		\begin{align*}
			\sum_{q\geq0}\alpha_{k,k-1}(T-s-qb_{k+1})
			&\leq\sum_{q\geq0}\alpha_{k,k-1}(T-qb_{k+1})
			=\sum_{q\geq0}\frac{1}{2^{E_k}}\sum_{j=0}^{E_k}\binom{E_k}{j}f_{k-1}(T-j-qb_{k+1})\\
			&=\frac{1}{2^{E_k}}\sum_{j=0}^{E_k}\binom{E_k}{j}\sum_{q\geq0}f_{k-1}(T-j-qb_{k+1})
			=\frac{1}{2^{E_k}}\sum_{j=0}^{E_k}\binom{E_k}{j}S_{k-1,\,b_{k+1}}(T-j).
		\end{align*}
		From~\eqref{eq:upper-step1} we get
		\begin{equation}\label{eq:upper-step2}
			D_{k+1}(n)\leq\frac{1}{b_2\dotsm b_k}\cdot\frac{1}{2^{E_k}}\sum_{j=0}^{E_k}\binom{E_k}{j}S_{k-1,\,b_{k+1}}(T-j).
		\end{equation}

		\item[$\scriptstyle\blacklozenge$] \emph{Step~3: the sums along the progression.} We apply Lemma~\ref{lem:sum-upper} with
		$r=k-1\geq1$ and $w=b_{k+1}$. For this $w$, by~\eqref{eq:UV} we have
		$U=\floor*{\frac{b_{k+1}}{2}}$ and, by~\eqref{eq:M-rec},
		$V=\floor*{\frac{b_{k+1}-1}{2}}=U-\varepsilon_{k+1}$. Let
		\begin{equation*}
			X
			\coloneqq T+U
			=n+M_k-(k-1)+\floor*{\frac{b_{k+1}}{2}}
			=n+\Bigl(M_k+\floor*{\frac{b_{k+1}}{2}}+1\Bigr)-k
			=n+M_{k+1}-k,
		\end{equation*}
		by~\eqref{eq:M-rec}. For $0\leq j\leq E_k$, Lemma~\ref{lem:sum-upper} at the point $T-j$ gives
		\begin{equation*}
			S_{k-1,\,b_{k+1}}(T-j)\leq\frac{1}{2b_{k+1}}\Bigl(f_k(X-j)+f_k(X-\varepsilon_{k+1}-j)\Bigr).
		\end{equation*}
		Since the weights $\binom{E_k}{j}$ are non-negative, we may insert this estimate
		in~\eqref{eq:upper-step2}; then, by~\eqref{eq:average-def} with $r=k$ and using
		$b_2\dotsm b_k\cdot b_{k+1}=b_2\dotsm b_{k+1}$,
		\begin{align}
			D_{k+1}(n)
			&\leq\frac{1}{b_2\dotsm b_k}\cdot\frac{1}{2^{E_k}}\sum_{j=0}^{E_k}\binom{E_k}{j}\frac{1}{2b_{k+1}}\Bigl(f_k(X-j)+f_k(X-\varepsilon_{k+1}-j)\Bigr)\notag\\
			&=\frac{1}{b_2\dotsm b_{k+1}}\cdot\frac12\Biggl(\frac{1}{2^{E_k}}\sum_{j=0}^{E_k}\binom{E_k}{j}f_k(X-j)
			+\frac{1}{2^{E_k}}\sum_{j=0}^{E_k}\binom{E_k}{j}f_k(X-\varepsilon_{k+1}-j)\Biggr)\notag\\
			&=\frac{1}{b_2\dotsm b_{k+1}}\cdot\frac12\Bigl(\alpha_{k,k}(X)+\alpha_{k,k}(X-\varepsilon_{k+1})\Bigr).
			\label{eq:upper-step3}
		\end{align}

		\item[$\scriptstyle\blacklozenge$] \emph{Step~4: collecting the terms.} By~\eqref{eq:upper-step3}, then by
		Lemma~\ref{lem:average}~(i) with $r=k$, and finally since $X=n+M_{k+1}-k$,
		\begin{equation*}
			D_{k+1}(n)
			\leq\frac{1}{b_2\dotsm b_{k+1}}\cdot\frac12\Bigl(\alpha_{k,k}(X)+\alpha_{k,k}(X-\varepsilon_{k+1})\Bigr)
			=\frac{1}{b_2\dotsm b_{k+1}}\,\alpha_{k+1,k}(X)
			=\frac{1}{b_2\dotsm b_{k+1}}\,\alpha_{k+1,k}\bigl(n+M_{k+1}-k\bigr),
		\end{equation*}
		which is~\eqref{eq:upper-goal}.
	\end{itemize}
\end{proof}

\subsection{Comparison with the upper bound of Agnarsson}\label{s:upper-agnarsson}
Both the upper bounds of Theorem~\ref{thm:main-upper} and that of
Agnarsson in~\eqref{eq:agnarsson-ours} have the factor $\frac{1}{b_2\dotsm b_k}$, so in the
comparison we leave it out.

\begin{prop}[Agnarsson's upper bound and Theorem~\ref{thm:main-upper}]\label{prop:agnarsson-upper}
	Let $k\in\N_{\geq2}$. For every $n\in\N$, the upper bounds of Theorem~\ref{thm:main-upper} are
	less than or equal to the upper bound of Agnarsson. The bound~\eqref{eq:upper-sum} is less if
	$k\geq3$ and $b_i\geq2$ for some $3\leq i\leq k$.
\end{prop}
\begin{proof}
	Let $n\in\N$. By Theorem~\ref{thm:main-upper}, the bound~\eqref{eq:upper-sum} is equal
	to~\eqref{eq:upper-average} and less than or equal to~\eqref{eq:upper-single}. So we have to
	prove that
	\begin{equation*}
		\binom{n+M_k}{k-1}\leq\binom{n+A_k}{k-1},
	\end{equation*}
	and that~\eqref{eq:upper-sum}, without its factor, is less than $\binom{n+A_k}{k-1}$ if
	$k\geq3$ and $b_i\geq2$ for some $3\leq i\leq k$.

	For every $3\leq i\leq k$ we have $\floor*{\frac{b_i}{2}}+1\leq b_i$, with equality only if
	$b_i\in\{1,2\}$. Indeed, if $b_i$ is even, then $\floor*{\frac{b_i}{2}}+1=\frac{b_i}{2}+1$,
	which is at most $b_i$ because $b_i\geq2$, with equality only for $b_i=2$; if $b_i$ is odd,
	then $\floor*{\frac{b_i}{2}}+1=\frac{b_i+1}{2}$, which is at most $b_i$ because $b_i\geq1$,
	with equality only for $b_i=1$. Hence, by~\eqref{eq:M-def} and~\eqref{eq:AB-ours},
	\begin{equation*}
		M_k
		=b_2+\sum_{i=3}^{k}\Bigl(\floor*{\frac{b_i}{2}}+1\Bigr)
		\leq b_2+\sum_{i=3}^{k}b_i
		=A_k,
	\end{equation*}
	with strict inequality if $b_i\geq3$ for some $3\leq i\leq k$. Moreover, by Pascal's rule,
	$\binom{y+1}{k-1}-\binom{y}{k-1}=\binom{y}{k-2}$ for every $y\in\N$; this is non-negative, and
	positive for $y\geq k-2$. Since $n+M_k\geq k-1$ by~\eqref{eq:M-lower}, going from
	$y=n+M_k$ up to $y=n+A_k$ we get
	\begin{equation*}
		\binom{n+M_k}{k-1}\leq\binom{n+A_k}{k-1},
	\end{equation*}
	with strict inequality if $M_k<A_k$.

	Let now $k\geq3$ and $b_i\geq2$ for some $3\leq i\leq k$. If $b_i\geq3$ for some such $i$, then
	$M_k<A_k$, and~\eqref{eq:upper-single} is less than Agnarsson's bound; hence so
	is~\eqref{eq:upper-sum}. Otherwise $b_i=2$ for some such $i$; then $b_i$ is even, and by
	Theorem~\ref{thm:main-upper} the bound~\eqref{eq:upper-sum} is less
	than~\eqref{eq:upper-single}, which is at most Agnarsson's bound.
\end{proof}

\begin{rem}[The terms of lower degree]\label{rem:upper-lower-degree}
	The terms of lower degree of~\eqref{eq:upper-sum} are much smaller in absolute value than
	those of the lower bound: the numbers $\frac{b_i-1}{2}$ of Theorem~\ref{thm:main-lower} are
	replaced by the numbers $-\frac{\varepsilon_i}{2}\in\{0,-\frac12\}$.
\end{rem}

\subsection{Proof of Theorem~\ref{thm:main}}
We have $\frac{d_k}{a_1\dotsm a_k}=\frac{1}{b_2\dotsm b_k}$ by Lemma~\ref{lem:basic}~(ii). Hence
the lower bound of~\eqref{eq:intro-main} is Theorem~\ref{thm:main-lower}, and the upper bound
is~\eqref{eq:upper-single} in Theorem~\ref{thm:main-upper}.

\section{Further comparisons}\label{s:compare}
This section compares the bounds of Theorems~\ref{thm:main-lower} and~\ref{thm:main-upper}
with the polynomial bounds of the first paper, and with bounds in terms of the least common
multiple of $a_1,\dotsc,a_k$.

\subsection{The polynomial lower bound}
Besides~\eqref{eq:intro-first}, the first paper contains a sharper lower bound, with terms of
lower degree (Theorem~3.5 there): for $k\geq2$ and every multiple $n\geq\Delta_k$ of $d_k$,
\begin{equation}\label{eq:previous-lower}
	D_k(n)>\frac{d_k}{a_1\dotsm a_k}\sum_{i=0}^{k-2}
	\bbinom{k-2}{i}\frac{(n-\Delta_k)^{k-1-i}}{(k-1-i)!},
\end{equation}
where $\bbinom{m}{i}$ is the symbol of Definition~\ref{dfn:coefficients} for the sequence
$c_l=b_l$; this is the notation of Section~2.4 of the first paper. We put the lower bound of
Theorem~\ref{thm:main-lower} next to it, both written with the factor
$\frac{1}{b_2\dotsm b_k}$ by Lemma~\ref{lem:basic}~(ii):
\begin{equation}\label{eq:previous-lower-ours}
	\begin{aligned}
		D_k(n)&\geq\frac{1}{b_2\dotsm b_k}\sum_{i=0}^{k-2}\bbinom[b-1]{k-2}{i}\binom{n-B_k-i}{k-1-i}
		&&\text{(Theorem~\ref{thm:main-lower})},\\
		D_k(n)&>\frac{1}{b_2\dotsm b_k}\sum_{i=0}^{k-2}\bbinom{k-2}{i}\frac{(n-\Delta_k)^{k-1-i}}{(k-1-i)!}
		&&\text{(the first paper)}.
	\end{aligned}
\end{equation}
Both hold for the same integers $n$, the multiples $n\geq\Delta_k$ of $d_k$, since
$\Delta_k=B_k+k-1$ by~\eqref{eq:AB-ours}. Comparing them term by term, the power
$\frac{(n-\Delta_k)^{k-1-i}}{(k-1-i)!}$ is replaced by the binomial coefficient
$\binom{n-B_k-i}{k-1-i}$, which has the same leading term, and the coefficient
$\bbinom{k-2}{i}$ is replaced by the coefficient $\bbinom[b-1]{k-2}{i}$, which is not larger. So it is
not clear a priori which of the two bounds is sharper.
Theorem~\ref{thm:compare-lower} shows that the bound of Theorem~\ref{thm:main-lower} always
is. For its proof we name the two sums of~\eqref{eq:previous-lower-ours} as functions of
$N=n-\Delta_k$. For $k\geq2$ and $N\in\N$ let
\begin{equation}\label{eq:GH-def}
	G_k(N)\coloneqq\sum_{i=0}^{k-2}\bbinom[b-1]{k-2}{i}\binom{N+k-1-i}{k-1-i},
	\qquad
	H_k(N)\coloneqq\sum_{i=0}^{k-2}\bbinom{k-2}{i}\frac{N^{k-1-i}}{(k-1-i)!}.
\end{equation}
Indeed, for $N=n-\Delta_k$ we have $N+k-1-i=n-B_k-i$ by~\eqref{eq:AB-ours}, so the two sums
of~\eqref{eq:previous-lower-ours} are $G_k(n-\Delta_k)$ and $H_k(n-\Delta_k)$. Each of the two
functions satisfies a recursion in $k$.

\begin{lem}\label{lem:recurrences}
	For every $N\in\N$, $G_2(N)=N+1$, $H_2(N)=N$ and, for $k\geq2$,
	\begin{align*}
		G_{k+1}(N)&=\sum_{t=0}^{N}G_k(t)+\frac{b_{k+1}-1}{2}\,G_k(N),\\
		H_{k+1}(N)&=\sum_{i=0}^{k-2}\bbinom{k-2}{i}\frac{N^{k-i}}{(k-i)!}+\frac{b_{k+1}}{2}\,H_k(N).
	\end{align*}
\end{lem}
\begin{proof}
	\emph{The case $k=2$.} Both sums have the single term $i=0$, with coefficient
	$\bbinom[b-1]{0}{0}=\bbinom{0}{0}=1$; hence $G_2(N)=\binom{N+1}{1}=N+1$ and $H_2(N)=\frac{N^1}{1!}=N$.

	\emph{The recurrence for $G$.} Let $k\geq2$. By Lemma~\ref{lem:pascal-sums} with $c=b-1$, $m=k-1$
	and $y_i\coloneqq\binom{N+k-i}{k-i}$ for $0\leq i\leq k-1$ (so that
	$y_{i+1}=\binom{N+k-1-i}{k-1-i}$),
	\begin{equation*}
		G_{k+1}(N)
		=\sum_{i=0}^{k-1}\bbinom[b-1]{k-1}{i}y_i
		=\sum_{i=0}^{k-2}\bbinom[b-1]{k-2}{i}\binom{N+k-i}{k-i}
		+\frac{b_{k+1}-1}{2}\sum_{i=0}^{k-2}\bbinom[b-1]{k-2}{i}\binom{N+k-1-i}{k-1-i}.
	\end{equation*}
	The second sum is $G_k(N)$, by definition. In the first sum, for each $0\leq i\leq k-2$,
	Lemma~\ref{lem:trunc}~(iii) with $k-i$ in place of $r$, $c=0$ and $N$ in place of $T$ gives
	\begin{equation*}
		\binom{N+k-i}{k-i}
		=f_{k-i}(N)
		=\sum_{t=0}^{N}f_{k-1-i}(t)
		=\sum_{t=0}^{N}\binom{t+k-1-i}{k-1-i},
	\end{equation*}
	since $N\geq0$ and $t\geq0$. Exchanging the two finite sums, we get
	\begin{equation*}
		\sum_{i=0}^{k-2}\bbinom[b-1]{k-2}{i}\binom{N+k-i}{k-i}
		=\sum_{t=0}^{N}\sum_{i=0}^{k-2}\bbinom[b-1]{k-2}{i}\binom{t+k-1-i}{k-1-i}
		=\sum_{t=0}^{N}G_k(t).
	\end{equation*}

	\emph{The recurrence for $H$.} By Lemma~\ref{lem:pascal-sums} with $c=b$, $m=k-1$ and
	$y_i\coloneqq\frac{N^{k-i}}{(k-i)!}$,
	\begin{equation*}
		H_{k+1}(N)=\sum_{i=0}^{k-2}\bbinom{k-2}{i}\frac{N^{k-i}}{(k-i)!}
		+\frac{b_{k+1}}{2}\sum_{i=0}^{k-2}\bbinom{k-2}{i}\frac{N^{k-1-i}}{(k-1-i)!},
	\end{equation*}
	and the second sum is $H_k(N)$, by definition.
\end{proof}

We also need the lower bound in Lemma~2.5 of~\cite{ballesio-arxiv}: for $j\in\N_{\geq1}$ and
$x\in\R_{\geq0}$, $\sum_{q=0}^{\floor{x}}(x-q)^j\geq\frac{x^{j+1}}{j+1}+\frac{x^j}{2}$. For an integer
$x=N$ the left-hand side is $\sum_{t=0}^{N}t^j$, so, dividing by $j!$,
\begin{equation}\label{eq:powers}
	\sum_{t=0}^{N}\frac{t^j}{j!}\geq\frac{N^{j+1}}{(j+1)!}+\frac12\cdot\frac{N^j}{j!}
	\qquad(j\in\N_{\geq1},\ N\in\N).
\end{equation}
\begin{thm}[The lower bound of the first paper and Theorem~\ref{thm:main-lower}]\label{thm:compare-lower}
	Let $k\in\N_{\geq2}$. For every integer $n\geq\Delta_k$, the lower bound of
	Theorem~\ref{thm:main-lower} is greater than the polynomial lower
	bound~\eqref{eq:previous-lower} of the first paper.
\end{thm}
\begin{proof}
	Let $n\geq\Delta_k$. By~\eqref{eq:previous-lower-ours} and~\eqref{eq:GH-def}, we have to prove
	that $G_k(n-\Delta_k)>H_k(n-\Delta_k)$. We prove that
	\begin{equation}\label{eq:GH-claim}
		G_k(N)>H_k(N)\qquad\text{for every }N\in\N
	\end{equation}
	by induction on $k$.

	\emph{Base case $k=2$.} By Lemma~\ref{lem:recurrences}, $G_2(N)=N+1>N=H_2(N)$.

	\emph{Inductive step.} Let $k\geq2$, and assume that
	\begin{equation}\label{eq:GH-hyp}
		G_k(t)>H_k(t)\qquad\text{for every }t\in\N.
	\end{equation}
	We have to prove that $G_{k+1}(N)>H_{k+1}(N)$ for every $N\in\N$. Let $N\in\N$.

	\begin{itemize}
		\item[$\scriptstyle\blacklozenge$] \emph{Step~1: the recurrence for $G$.} By Lemma~\ref{lem:recurrences}, and then by the
		inductive hypothesis~\eqref{eq:GH-hyp} for $t=0,1,\dotsc,N$, where
		$\frac{b_{k+1}-1}{2}\geq0$ since $b_{k+1}\geq1$,
		\begin{equation}\label{eq:GH-step1}
			G_{k+1}(N)
			=\sum_{t=0}^{N}G_k(t)+\frac{b_{k+1}-1}{2}\,G_k(N)
			>\sum_{t=0}^{N}H_k(t)+\frac{b_{k+1}-1}{2}\,H_k(N).
		\end{equation}
		The inequality is strict because the sum contains at least the term $t=0$.

		\item[$\scriptstyle\blacklozenge$] \emph{Step~2: the sums of powers.} By~\eqref{eq:GH-def} and exchanging the two finite sums,
		\begin{equation*}
			\sum_{t=0}^{N}H_k(t)
			=\sum_{t=0}^{N}\sum_{i=0}^{k-2}\bbinom{k-2}{i}\frac{t^{k-1-i}}{(k-1-i)!}
			=\sum_{i=0}^{k-2}\bbinom{k-2}{i}\sum_{t=0}^{N}\frac{t^{k-1-i}}{(k-1-i)!}.
		\end{equation*}
		For $0\leq i\leq k-2$ the exponent $j=k-1-i$ satisfies $j\geq1$, so~\eqref{eq:powers}
		applies to the inner sum:
		\begin{equation*}
			\sum_{t=0}^{N}\frac{t^{k-1-i}}{(k-1-i)!}\geq\frac{N^{k-i}}{(k-i)!}+\frac12\cdot\frac{N^{k-1-i}}{(k-1-i)!}.
		\end{equation*}
		The coefficients $\bbinom{k-2}{i}$ are non-negative (Lemma~\ref{lem:sym-poly}~(ii)), so we may
		insert this estimate in the previous identity; then, by~\eqref{eq:GH-def},
		\begin{equation}\label{eq:GH-step2}
			\sum_{t=0}^{N}H_k(t)
			\geq\sum_{i=0}^{k-2}\bbinom{k-2}{i}\frac{N^{k-i}}{(k-i)!}
			+\frac12\sum_{i=0}^{k-2}\bbinom{k-2}{i}\frac{N^{k-1-i}}{(k-1-i)!}
			=\sum_{i=0}^{k-2}\bbinom{k-2}{i}\frac{N^{k-i}}{(k-i)!}+\frac12\,H_k(N).
		\end{equation}

		\item[$\scriptstyle\blacklozenge$] \emph{Step~3: the recurrence for $H$.} By~\eqref{eq:GH-step1}, then by~\eqref{eq:GH-step2},
		then since $\frac12+\frac{b_{k+1}-1}{2}=\frac{b_{k+1}}{2}$, and finally by
		Lemma~\ref{lem:recurrences},
		\begin{align*}
			G_{k+1}(N)
			&>\sum_{t=0}^{N}H_k(t)+\frac{b_{k+1}-1}{2}\,H_k(N)
			\geq\sum_{i=0}^{k-2}\bbinom{k-2}{i}\frac{N^{k-i}}{(k-i)!}+\frac12\,H_k(N)+\frac{b_{k+1}-1}{2}\,H_k(N)\\
			&=\sum_{i=0}^{k-2}\bbinom{k-2}{i}\frac{N^{k-i}}{(k-i)!}+\frac{b_{k+1}}{2}\,H_k(N)
			=H_{k+1}(N),
		\end{align*}
		which is what we had to prove.
	\end{itemize}
\end{proof}

\subsection{The polynomial upper bound}
We put the upper bound~\eqref{eq:upper-single} of Theorem~\ref{thm:main-upper} next to the upper
bound of the first paper (Theorem~3.2 there), written with the factor $\frac{1}{b_2\dotsm b_k}$:
for every multiple $n$ of $d_k$,
\begin{equation}\label{eq:previous-upper}
	\begin{aligned}
		D_k(n)&\leq\frac{1}{b_2\dotsm b_k}\binom{n+M_k}{k-1}
		&&\text{(Theorem~\ref{thm:main-upper})},\\
		D_k(n)&\leq\frac{1}{b_2\dotsm b_k}\cdot\frac{(n+\nabla_k)^{k-1}}{(k-1)!}
		&&\text{(the first paper)}.
	\end{aligned}
\end{equation}
Unlike the two lower bounds, neither of these upper bounds is always the sharper one: the
outcome depends on $k$ and on the parities of $b_3,\dotsc,b_k$. There are three cases.
\begin{itemize}
	\item If $b_3,\dotsc,b_k$ are all odd, the bound of Theorem~\ref{thm:main-upper} is never
	      weaker: this is Proposition~\ref{prop:compare-upper}~(i).
	\item If $k=3$ and $b_3$ is even, the two bounds are equal: this is
	      Proposition~\ref{prop:compare-upper}~(ii).
	\item If $k\geq4$ and at least one of $b_3,\dotsc,b_k$ is even, either bound can be the
	      sharper one: see the examples of Remark~\ref{rem:compare-upper}.
\end{itemize}

\begin{prop}[The upper bound of the first paper and Theorem~\ref{thm:main-upper}]\label{prop:compare-upper}
	Let $k\in\N_{\geq2}$ and $n\in\N$.
	\begin{enumerate}
		\item[\textup{(i)}] If $b_3,\dotsc,b_k$ are all odd, then the upper bounds of
		      Theorem~\ref{thm:main-upper} are less than or equal to the upper bound of the first
		      paper.
		\item[\textup{(ii)}] If $k=3$ and $b_3$ is even, then the upper bound~\eqref{eq:upper-sum}
		      of Theorem~\ref{thm:main-upper} is equal to the upper bound of the first paper.
	\end{enumerate}
\end{prop}
\begin{proof}
	In both parts we leave out the common factor $\frac{1}{b_2\dotsm b_k}$.
	\begin{enumerate}
		\item[\textup{(i)}] By Theorem~\ref{thm:main-upper}, if $b_3,\dotsc,b_k$ are all odd the
		      three bounds~\eqref{eq:upper-sum}, \eqref{eq:upper-average}
		      and~\eqref{eq:upper-single} are equal. So we have to prove that
		      \begin{equation*}
			      \binom{n+M_k}{k-1}\leq\frac{(n+\nabla_k)^{k-1}}{(k-1)!}.
		      \end{equation*}
		      For an odd $b_i$ we have $\floor*{\frac{b_i}{2}}+1=\frac{b_i}{2}+\frac12$, hence
		      by~\eqref{eq:M-def}
		      \begin{equation*}
			      M_k
			      =b_2+\sum_{i=3}^{k}\Bigl(\frac{b_i}{2}+\frac12\Bigr)
			      =\nabla_k+\frac{k-2}{2}.
		      \end{equation*}
		      The numerator of $\binom{n+M_k}{k-1}$ is the product of the $k-1$ integers $n+M_k-j$
		      with $0\leq j\leq k-2$, which are positive since $M_k\geq k-1$
		      by~\eqref{eq:M-lower}. By the inequality between the geometric and the arithmetic
		      mean,
		      \begin{equation*}
			      \prod_{j=0}^{k-2}(n+M_k-j)
			      \leq\Biggl(\frac{1}{k-1}\sum_{j=0}^{k-2}(n+M_k-j)\Biggr)^{k-1}.
		      \end{equation*}
		      The sum has $k-1$ terms, and $\sum_{j=0}^{k-2}j=\frac{(k-2)(k-1)}{2}$; hence, and
		      then by the value of $M_k$ found above, the arithmetic mean is
		      \begin{equation*}
			      \frac{1}{k-1}\sum_{j=0}^{k-2}(n+M_k-j)
			      =\frac{1}{k-1}\Bigl((k-1)(n+M_k)-\frac{(k-2)(k-1)}{2}\Bigr)
			      =n+M_k-\frac{k-2}{2}
			      =n+\nabla_k.
		      \end{equation*}
		      So the product is at most $(n+\nabla_k)^{k-1}$.
		      Dividing by $(k-1)!$ we get the inequality.
		\item[\textup{(ii)}] By Theorem~\ref{thm:main-upper}, the bound~\eqref{eq:upper-sum} is
		      equal to~\eqref{eq:upper-average}. For $k=3$ and $b_3$ even we have $E_3=1$, so the
		      sum in~\eqref{eq:upper-average} is
		      $\frac12\bigl(\binom{n+M_3}{2}+\binom{n+M_3-1}{2}\bigr)$, and we have to prove that
		      \begin{equation*}
			      \frac12\Biggl(\binom{n+M_3}{2}+\binom{n+M_3-1}{2}\Biggr)=\frac{(n+\nabla_3)^{2}}{2!}.
		      \end{equation*}
		      By~\eqref{eq:M-def}, $M_3=b_2+\frac{b_3}{2}+1=\nabla_3+1$. Setting $y\coloneqq n+M_3$,
		      so that $y-1=n+\nabla_3$,
		      \begin{equation*}
			      \frac12\Biggl(\binom{y}{2}+\binom{y-1}{2}\Biggr)
			      =\frac{y(y-1)+(y-1)(y-2)}{4}
			      =\frac{(y-1)(2y-2)}{4}
			      =\frac{(y-1)^2}{2}
			      =\frac{(n+\nabla_3)^2}{2!}.\qedhere
		      \end{equation*}
	\end{enumerate}
\end{proof}

\begin{rem}\label{rem:compare-upper}
	\emph{What happens for $k\geq4$ when some $b_i$ is even?} The two upper bounds are not
	comparable. For instance, for $k=4$ and $n=0$, where $D_4(0)=1$:
	\begin{itemize}
		\item if $(a_1,a_2,a_3,a_4)=(1,1,2,2)$, so that $b_3=b_4=2$, the
		      bound~\eqref{eq:upper-sum} is $\frac{19}{16}$ and that of the first paper is
		      $\frac98=\frac{18}{16}$;
		\item if $(a_1,a_2,a_3,a_4)=(1,1,1,2)$, so that $b_3=1$ and $b_4=2$, the
		      bound~\eqref{eq:upper-sum} is $\frac54=\frac{120}{96}$ and that of the first paper
		      is $\frac{125}{96}$.
	\end{itemize}
\end{rem}

\subsection{Bounds in terms of the least common multiple}
Let $d_k=1$ and $\mu\coloneqq\lcm(a_1,\dotsc,a_k)$, and for $n\in\N$ write $n=h\mu+\rho$ with
$h\in\N$ and $0\leq\rho<\mu$. E.~Samsonadze proved the exact formula
(Theorem~3.1 and formula~(3.14) of~\cite{samsonadze-gmj})
\begin{equation}\label{eq:samsonadze}
	D_k(n)=\sum_{j=0}^{\sigma}\ell_j\binom{h+k-1-j}{k-1},
	\qquad
	\sum_{j=0}^{\sigma}\ell_j=\frac{\mu^{k-1}}{a_1\dotsm a_k},
\end{equation}
where $\sigma\coloneqq\floor*{k-\frac{a_1+\dotsb+a_k+\rho}{\mu}}\leq k-1$ and $\ell_j$ is the number
of solutions of $a_1x_1+\dotsb+a_kx_k=\rho+j\mu$ with $0\leq x_i<\frac{\mu}{a_i}$ for every $i$
(this description of $\ell_j$ is in the proof of Theorem~2.2 of~\cite{samsonadze-gmj}).
The numbers $\ell_j$ are non-negative, and for $0\leq j\leq\sigma$ we have
$\binom{h+k-1-\sigma}{k-1}\leq\binom{h+k-1-j}{k-1}\leq\binom{h+k-1}{k-1}$. Replacing each
binomial coefficient of~\eqref{eq:samsonadze} by the smallest and by the largest one, and then
using the value of $\sum_{j}\ell_j$, we get
\begin{equation}\label{eq:lcm-bounds}
	\frac{\mu^{k-1}}{a_1\dotsm a_k}\binom{h+k-1-\sigma}{k-1}
	\leq D_k(n)
	\leq\frac{\mu^{k-1}}{a_1\dotsm a_k}\binom{h+k-1}{k-1}.
\end{equation}
These bounds are of a different nature from those of Theorems~\ref{thm:main-lower}
and~\ref{thm:main-upper}: they involve the least common multiple $\mu$ in place of the numbers
$b_i$, and as functions of $n$ they are constant on the intervals between two consecutive
multiples of $\mu$, apart from the dependence of $\sigma$ on $\rho$.

\begin{rem}\label{rem:lcm}
	The bounds~\eqref{eq:lcm-bounds} and those of Theorems~\ref{thm:main-lower}
	and~\ref{thm:main-upper} are not comparable. For instance:
	\begin{itemize}
		\item for $(a_1,a_2,a_3)=(6,10,15)$ we have $b_2=b_3=\mu=30$, and for $n=60$, so that
		      $h=2$, $\rho=0$ and $\sigma=1$, \eqref{eq:lcm-bounds} gives $3\leq D_3(60)\leq6$,
		      while Theorems~\ref{thm:main-lower} and~\ref{thm:main-upper} give
		      $\frac{182}{225}\leq D_3(60)\leq\frac{49}{8}$; here $D_3(60)=6$;
		\item for $(a_1,a_2,a_3)=(3,5,7)$ we have $b_2=15$, $b_3=7$ and $\mu=105$, and for
		      $n=100$, so that $h=0$, $\rho=100$ and $\sigma=1$, \eqref{eq:lcm-bounds} gives
		      $0\leq D_3(100)\leq105$, while Theorems~\ref{thm:main-lower}
		      and~\ref{thm:main-upper} give $\frac{637}{15}\leq D_3(100)\leq\frac{1003}{15}$; here
		      $D_3(100)=55$.
	\end{itemize}
	Roughly speaking, \eqref{eq:lcm-bounds} is sharper when $\mu$ is not much larger than the
	numbers $b_i$, and the bounds of Theorems~\ref{thm:main-lower} and~\ref{thm:main-upper} are
	sharper when $\mu$ is large.
\end{rem}

\subsection{The case \texorpdfstring{$a_1=1$}{a1=1}}\label{s:a1}
If $a_1=1$, then $d_k=1$ and $b_i=a_i$ for every $i\geq2$, as shown in Section~5.1 of the first
paper. Hence
\begin{equation*}
	b_2\dotsm b_k=a_1\dotsm a_k,\qquad \Delta_k=0,\qquad B_k=-(k-1),\qquad
	A_k=\sum_{i=2}^{k}a_i,\qquad M_k=a_2+\sum_{i=3}^{k}\Bigl(\floor*{\frac{a_i}{2}}+1\Bigr).
\end{equation*}
and the sequence $b-1$ of~\eqref{eq:coeff-b-1} is the sequence $c_l=a_l-1$, which we denote by
$a-1$. Theorems~\ref{thm:main-lower} and~\ref{thm:main-upper} then hold for every $n\in\N$ and read
\begin{equation*}
	\frac{1}{a_1\dotsm a_k}\sum_{i=0}^{k-2}
	\bbinom[a-1]{k-2}{i}\binom{n+k-1-i}{k-1-i}
	\leq D_k(n)
	\leq\frac{1}{a_1\dotsm a_k}\binom{n+M_k}{k-1}.
\end{equation*}
The sharper upper bound~\eqref{eq:upper-sum} also applies, with $\varepsilon_i=1$ exactly when $a_i$
is even. The term $i=0$ of the lower bound and the upper bound with $A_k$ in place of $M_k$ are the
bounds of Agnarsson for $a_1=1$ (remark after Corollary~3.5 of~\cite{agnarsson}). By
Theorem~\ref{thm:compare-lower}, the lower bound is greater than the refined lower bound of
Blom and Fr\"oberg~\cite{blom-froberg}, which is~\eqref{eq:previous-lower} for $a_1=1$.

\subsection{The distance between the upper and the lower bound}\label{s:distance}
An upper and a lower bound for $D_k(n)$ locate it in an interval, whose length is the upper
bound minus the lower bound. Here we compare this length for the sharpest bounds of
the first paper and for the bounds of this paper. We denote the four bounds as follows, all
with the factor $\frac{1}{b_2\dotsm b_k}$:
\begin{align*}
	U_k^{\mathrm{I}}(n)&\coloneqq\frac{1}{b_2\dotsm b_k}\cdot\frac{(n+\nabla_k)^{k-1}}{(k-1)!}
	&&\text{(the first paper, see~\eqref{eq:previous-upper})},\\
	L_k^{\mathrm{I}}(n)&\coloneqq\frac{1}{b_2\dotsm b_k}\sum_{i=0}^{k-2}\bbinom{k-2}{i}\frac{(n-\Delta_k)^{k-1-i}}{(k-1-i)!}
	&&\text{(the first paper, see~\eqref{eq:previous-lower})},\\
	U_k^{\mathrm{II}}(n)&\coloneqq\frac{1}{b_2\dotsm b_k}\sum_{i=0}^{k-2}\bbinom[-\varepsilon]{k-2}{i}\binom{n+M_k-i}{k-1-i}
	&&\text{(this paper, see~\eqref{eq:upper-sum})},\\
	L_k^{\mathrm{II}}(n)&\coloneqq\frac{1}{b_2\dotsm b_k}\sum_{i=0}^{k-2}\bbinom[b-1]{k-2}{i}\binom{n-B_k-i}{k-1-i}
	&&\text{(this paper, Theorem~\ref{thm:main-lower})}.
\end{align*}
The first two are written with powers of $n$ and the other two with binomial coefficients, so
they cannot be compared at sight. But each of the four is a polynomial in $n$ of degree $k-1$,
and expanding them in powers of $n$ we can compute the two distances.

\begin{prop}[The distance between the upper and the lower bound]\label{prop:distance}
	Let $k\in\N_{\geq3}$. Then
	\begin{equation}\label{eq:distance}
		\begin{aligned}
			U_k^{\mathrm{I}}(n)-L_k^{\mathrm{I}}(n)
			&=\frac{b_2+\Delta_k}{b_2\dotsm b_k}\cdot\frac{n^{k-2}}{(k-2)!}+O\bigl(n^{k-3}\bigr)
			&&\text{(the first paper)},\\
			U_k^{\mathrm{II}}(n)-L_k^{\mathrm{II}}(n)
			&=\frac{b_2+\Delta_k-1}{b_2\dotsm b_k}\cdot\frac{n^{k-2}}{(k-2)!}+O\bigl(n^{k-3}\bigr)
			&&\text{(this paper)}.
		\end{aligned}
	\end{equation}
\end{prop}
\begin{proof}
	We have to prove~\eqref{eq:distance}. We first compute, for a power and for a binomial
	coefficient, the two terms of highest degree in $n$; then we apply this to each of the four
	bounds; finally we subtract.

	\emph{The first two terms of a power.} Let $r\in\N_{\geq1}$ and $z\in\R$. By the binomial
	theorem, $(n+z)^r=n^r+rz\,n^{r-1}+O\bigl(n^{r-2}\bigr)$. Dividing by $r!$, and using
	$\frac{r}{r!}=\frac{1}{(r-1)!}$ in the second term,
	\begin{equation}\label{eq:power-two-terms}
		\frac{(n+z)^r}{r!}
		=\frac{n^r}{r!}+z\,\frac{n^{r-1}}{(r-1)!}+O\bigl(n^{r-2}\bigr),
	\end{equation}
	where for $r=1$ there is no $O$-term.

	\emph{The first two terms of a binomial coefficient.} Let $r\in\N_{\geq1}$ and $z\in\R$. By
	definition, $\binom{n+z}{r}$ is the product of the $r$ factors $n+z-j$, with
	$0\leq j\leq r-1$, divided by $r!$. We expand the product, choosing in each factor either $n$
	or the constant $z-j$:
	\begin{itemize}
		\item choosing $n$ in all the $r$ factors gives $n^r$;
		\item choosing the constant $z-j$ in one factor and $n$ in the other $r-1$ factors gives
		      $(z-j)\,n^{r-1}$, once for each $0\leq j\leq r-1$;
		\item choosing the constant in two or more factors gives terms of degree at most $r-2$.
	\end{itemize}
	Hence, and then since $\sum_{j=0}^{r-1}j=\frac{r(r-1)}{2}$,
	\begin{align*}
		(n+z)(n+z-1)\dotsm(n+z-r+1)
		&=n^r+\bigl(z+(z-1)+\dotsb+(z-r+1)\bigr)\,n^{r-1}+O\bigl(n^{r-2}\bigr)\\
		&=n^r+\Bigl(rz-\frac{r(r-1)}{2}\Bigr)n^{r-1}+O\bigl(n^{r-2}\bigr).
	\end{align*}
	Dividing by $r!$, and using $\frac{r}{r!}=\frac{1}{(r-1)!}$ in the second term,
	\begin{equation}\label{eq:binom-two-terms}
		\binom{n+z}{r}
		=\frac{n^r}{r!}+\Bigl(z-\frac{r-1}{2}\Bigr)\frac{n^{r-1}}{(r-1)!}+O\bigl(n^{r-2}\bigr),
	\end{equation}
	where for $r=1$ there is no $O$-term.

	\emph{The upper bound of the first paper.} By~\eqref{eq:power-two-terms} with $r=k-1$ and
	$z=\nabla_k$,
	\begin{equation}\label{eq:expansion-UI}
		U_k^{\mathrm{I}}(n)
		=\frac{1}{b_2\dotsm b_k}\cdot\frac{(n+\nabla_k)^{k-1}}{(k-1)!}
		=\frac{1}{b_2\dotsm b_k}\biggl(\frac{n^{k-1}}{(k-1)!}+\nabla_k\,\frac{n^{k-2}}{(k-2)!}\biggr)+O\bigl(n^{k-3}\bigr).
	\end{equation}

	\emph{The lower bound of the first paper.} Apart from the factor $\frac{1}{b_2\dotsm b_k}$,
	$L_k^{\mathrm{I}}(n)$ is a sum over $0\leq i\leq k-2$, whose term with index $i$ is a
	polynomial in $n$ of degree $k-1-i$. By~\eqref{eq:power-two-terms} with $z=-\Delta_k$, with
	$r=k-1$ for the term $i=0$ and with $r=k-2$ for the term $i=1$, and since
	$\bbinom{k-2}{0}=1$ and $\bbinom{k-2}{1}=\sum_{l=3}^{k}\frac{b_l}{2}=\nabla_k-b_2$
	(Definition~\ref{dfn:constants}),
	\begin{align*}
		\bbinom{k-2}{0}\frac{(n-\Delta_k)^{k-1}}{(k-1)!}
		&=\frac{n^{k-1}}{(k-1)!}-\Delta_k\,\frac{n^{k-2}}{(k-2)!}+O\bigl(n^{k-3}\bigr)
		&&\text{(the term $i=0$)},\\
		\bbinom{k-2}{1}\frac{(n-\Delta_k)^{k-2}}{(k-2)!}
		&=(\nabla_k-b_2)\,\frac{n^{k-2}}{(k-2)!}+O\bigl(n^{k-3}\bigr)
		&&\text{(the term $i=1$)}.
	\end{align*}
	The terms with $i\geq2$ have degree $k-1-i\leq k-3$, so they are $O\bigl(n^{k-3}\bigr)$.
	Adding all the terms,
	\begin{equation}\label{eq:expansion-LI}
		L_k^{\mathrm{I}}(n)
		=\frac{1}{b_2\dotsm b_k}\biggl(\frac{n^{k-1}}{(k-1)!}+\bigl(\nabla_k-(b_2+\Delta_k)\bigr)\frac{n^{k-2}}{(k-2)!}\biggr)+O\bigl(n^{k-3}\bigr).
	\end{equation}

	\emph{The upper bound of this paper.} Apart from the factor $\frac{1}{b_2\dotsm b_k}$,
	$U_k^{\mathrm{II}}(n)$ is a sum over $0\leq i\leq k-2$, whose term with index $i$ is a
	polynomial in $n$ of degree $k-1-i$. By~\eqref{eq:binom-two-terms}, with $r=k-1$ and $z=M_k$
	for the term $i=0$ and with $r=k-2$ and $z=M_k-1$ for the term $i=1$, and since
	$\bbinom[-\varepsilon]{k-2}{0}=1$ and $\bbinom[-\varepsilon]{k-2}{1}=-\frac{E_k}{2}$
	by~\eqref{eq:coeff-eps},
	\begin{align*}
		\bbinom[-\varepsilon]{k-2}{0}\binom{n+M_k}{k-1}
		&=\frac{n^{k-1}}{(k-1)!}+\Bigl(M_k-\frac{k-2}{2}\Bigr)\frac{n^{k-2}}{(k-2)!}+O\bigl(n^{k-3}\bigr)
		&&\text{(the term $i=0$)},\\
		\bbinom[-\varepsilon]{k-2}{1}\binom{n+M_k-1}{k-2}
		&=-\frac{E_k}{2}\cdot\frac{n^{k-2}}{(k-2)!}+O\bigl(n^{k-3}\bigr)
		&&\text{(the term $i=1$)}.
	\end{align*}
	The terms with $i\geq2$ have degree $k-1-i\leq k-3$, so they are $O\bigl(n^{k-3}\bigr)$. To
	add the two coefficients of $\frac{n^{k-2}}{(k-2)!}$ we write $M_k$ in terms of $\nabla_k$.
	By~\eqref{eq:eps-E-def}, $\floor*{\frac{b_l}{2}}=\frac{b_l-1+\varepsilon_l}{2}$ for every $l$;
	hence, by~\eqref{eq:M-def},
	\begin{equation*}
		M_k
		=b_2+\sum_{l=3}^{k}\Bigl(\frac{b_l}{2}+\frac{1+\varepsilon_l}{2}\Bigr)
		=\nabla_k+\frac{k-2}{2}+\frac{E_k}{2},
		\qquad\text{so that}\qquad
		M_k-\frac{k-2}{2}-\frac{E_k}{2}=\nabla_k .
	\end{equation*}
	Adding all the terms,
	\begin{equation}\label{eq:expansion-UII}
		U_k^{\mathrm{II}}(n)
		=\frac{1}{b_2\dotsm b_k}\biggl(\frac{n^{k-1}}{(k-1)!}+\nabla_k\,\frac{n^{k-2}}{(k-2)!}\biggr)+O\bigl(n^{k-3}\bigr).
	\end{equation}

	\emph{The lower bound of this paper.} Apart from the factor $\frac{1}{b_2\dotsm b_k}$,
	$L_k^{\mathrm{II}}(n)$ is a sum over $0\leq i\leq k-2$, whose term with index $i$ is a
	polynomial in $n$ of degree $k-1-i$. By~\eqref{eq:binom-two-terms}, with $r=k-1$ and
	$z=-B_k$ for the term $i=0$ and with $r=k-2$ and $z=-B_k-1$ for the term $i=1$, and since
	$\bbinom[b-1]{k-2}{0}=1$ and $\bbinom[b-1]{k-2}{1}=\sum_{l=3}^{k}\frac{b_l-1}{2}$
	by~\eqref{eq:coeff-b-1},
	\begin{align*}
		\bbinom[b-1]{k-2}{0}\binom{n-B_k}{k-1}
		&=\frac{n^{k-1}}{(k-1)!}+\Bigl(-B_k-\frac{k-2}{2}\Bigr)\frac{n^{k-2}}{(k-2)!}+O\bigl(n^{k-3}\bigr)
		&&\text{(the term $i=0$)},\\
		\bbinom[b-1]{k-2}{1}\binom{n-B_k-1}{k-2}
		&=\sum_{l=3}^{k}\frac{b_l-1}{2}\cdot\frac{n^{k-2}}{(k-2)!}+O\bigl(n^{k-3}\bigr)
		&&\text{(the term $i=1$)}.
	\end{align*}
	The terms with $i\geq2$ have degree $k-1-i\leq k-3$, so they are $O\bigl(n^{k-3}\bigr)$. The
	two coefficients of $\frac{n^{k-2}}{(k-2)!}$ add up, by~\eqref{eq:AB-ours} and then by
	Definition~\ref{dfn:constants}, to
	\begin{equation*}
		-B_k-\frac{k-2}{2}+\sum_{l=3}^{k}\frac{b_l-1}{2}
		=-\Delta_k+(k-1)-\frac{k-2}{2}+\sum_{l=3}^{k}\frac{b_l}{2}-\frac{k-2}{2}
		=\nabla_k-(b_2+\Delta_k-1).
	\end{equation*}
	Adding all the terms,
	\begin{equation}\label{eq:expansion-LII}
		L_k^{\mathrm{II}}(n)
		=\frac{1}{b_2\dotsm b_k}\biggl(\frac{n^{k-1}}{(k-1)!}+\bigl(\nabla_k-(b_2+\Delta_k-1)\bigr)\frac{n^{k-2}}{(k-2)!}\biggr)+O\bigl(n^{k-3}\bigr).
	\end{equation}

	\emph{The distances.} Subtracting~\eqref{eq:expansion-LI} from~\eqref{eq:expansion-UI}, and
	\eqref{eq:expansion-LII} from~\eqref{eq:expansion-UII}, the terms of degree $k-1$ cancel, and
	the coefficients of $\frac{n^{k-2}}{(k-2)!}$ are
	$\nabla_k-\bigl(\nabla_k-(b_2+\Delta_k)\bigr)=b_2+\Delta_k$ and
	$\nabla_k-\bigl(\nabla_k-(b_2+\Delta_k-1)\bigr)=b_2+\Delta_k-1$. This is~\eqref{eq:distance}.
\end{proof}

\begin{rem}\label{rem:distance}
	\begin{enumerate}
		\item[(a)] \emph{What changes from the first paper to this one?} In the term of degree
		      $k-2$ of the distance, the coefficient $b_2+\Delta_k$ becomes $b_2+\Delta_k-1$.
		      The two upper bounds have the same first two terms, by~\eqref{eq:expansion-UI}
		      and~\eqref{eq:expansion-UII},
		      so the whole gain comes from the lower bound; this agrees with Theorem~\ref{thm:compare-lower}.
		\item[(b)] \emph{The case $k=2$.} Here the four bounds are
		      $U_2^{\mathrm{I}}(n)=U_2^{\mathrm{II}}(n)=\frac{n+b_2}{b_2}$,
		      $L_2^{\mathrm{I}}(n)=\frac{n-\Delta_2}{b_2}$ and
		      $L_2^{\mathrm{II}}(n)=\frac{n-\Delta_2+1}{b_2}$, so the two distances are exactly
		      $\frac{b_2+\Delta_2}{b_2}$ and $\frac{b_2+\Delta_2-1}{b_2}$, in agreement
		      with~\eqref{eq:distance}.
		\item[(c)] \emph{When does the term of degree $k-2$ of the distance vanish?} In the
		      first paper, never, since $b_2\geq1$ and $\Delta_k\geq0$. In this paper, exactly
		      when $a_1=a_2=1$: indeed $b_2+\Delta_k-1=0$ if and only if $b_2=1$ and
		      $\Delta_k=0$; and $b_2=\lcm(a_1,a_2)=1$ if and only if $a_1=a_2=1$, in which case
		      $b_i=a_i$ for every $i$ (Section~\ref{s:a1}), so that $\Delta_k=0$.
	\end{enumerate}
\end{rem}

\section{Final remarks}\label{s:final}

\subsection{Refined bounds at the last step}
In the proofs of Theorems~\ref{thm:main-lower} and~\ref{thm:main-upper}, the integer $s$ of
Lemma~\ref{lem:recursion} is estimated by $0\leq s\leq b_{k+1}-a_{k+1}$, uniformly in $n$. As in
Section~6.1 of the first paper, keeping $s$ at the last step of the induction gives sharper
bounds, which depend on $n$ through $s$.

\begin{prop}\label{prop:refined}
	Let $k\in\N_{\geq2}$, let $n\in\N$ be a multiple of $d_k$, and let $s$ be the integer with
	$0\leq s<b_k$, $a_k\mid s$ and $d_{k-1}\mid n-s$, which exists and is unique by Lemma~2.2 of
	the first paper. For $k=2$ we read $\Delta_1$ as the empty sum $\sum_{i=2}^{1}(b_i-a_i)=0$.
	\begin{enumerate}
		\item[\textup{(i)}] If $n\geq s$, then
		      \begin{equation*}
			      D_k(n)
			      \leq\frac{1}{b_2\dotsm b_k}\sum_{i=0}^{k-2}\bbinom[-\varepsilon]{k-2}{i}\binom{n-s+M_k-i}{k-1-i}
			      \leq\frac{1}{b_2\dotsm b_k}\binom{n-s+M_k}{k-1}.
		      \end{equation*}
		\item[\textup{(ii)}] If $n\geq s+\Delta_{k-1}$, then
		      \begin{equation*}
			      D_k(n)\geq\frac{1}{b_2\dotsm b_k}\sum_{i=0}^{k-2}\bbinom[b-1]{k-2}{i}\binom{n-s-\Delta_{k-1}+k-1-i}{k-1-i}.
		      \end{equation*}
	\end{enumerate}
\end{prop}
In words: (i) is Theorem~\ref{thm:main-upper} with $n-s$ in place of $n$, and (ii) is
Theorem~\ref{thm:main-lower} with $s+\Delta_{k-1}$ in place of $\Delta_k$, since
$\binom{n-B_k-i}{k-1-i}=\binom{n-\Delta_k+k-1-i}{k-1-i}$ by~\eqref{eq:AB-ours}.
\begin{proof}
	\emph{The case $k=2$.} Here $M_2=b_2$, $E_2=0$ and $\Delta_1=0$, the sums have the single
	term $i=0$ with coefficient $1$, and both (i) and (ii) assume $n\geq s$. So we have to prove
	that
	\begin{equation*}
		\frac{n-s+1}{b_2}\leq D_2(n)\leq\frac{n-s+b_2}{b_2}.
	\end{equation*}
	Since $n\geq s$, Proposition~3.1 of the first paper gives
	$D_2(n)=\floor*{\frac{n-s}{b_2}}+1$. For every $y\in\R$ we have $y<\floor{y}+1\leq y+1$; with
	$y=\frac{n-s}{b_2}$ this gives $\frac{n-s}{b_2}<D_2(n)\leq\frac{n-s+b_2}{b_2}$. The first
	inequality reads $n-s<b_2D_2(n)$; both sides are integers, hence $n-s+1\leq b_2D_2(n)$.

	\emph{The case $k\geq3$.} By Lemma~\ref{lem:recursion}, for the passage from $k-1$ to $k$,
	\begin{equation}\label{eq:refined-rec}
		D_k(n)=\sum_{q=0}^{\floor*{\frac{n-s}{b_k}}}D_{k-1}(n-s-qb_k),
	\end{equation}
	where every argument $n-s-qb_k$ is a non-negative multiple of $d_{k-1}$. By
	Theorems~\ref{thm:main-lower} and~\ref{thm:main-upper} for $D_{k-1}$, in the
	forms~\eqref{eq:claim-lower} and~\eqref{eq:claim-upper}, we know that, for every multiple
	$m\in\N$ of $d_{k-1}$,
	\begin{equation}\label{eq:refined-known}
		\frac{1}{b_2\dotsm b_{k-1}}\sum_{i=0}^{k-3}\bbinom[b-1]{k-3}{i}f_{k-2-i}(m-\Delta_{k-1})
		\leq D_{k-1}(m)
		\leq\frac{1}{b_2\dotsm b_{k-1}}\,\alpha_{k-1,k-2}\bigl(m+M_{k-1}-(k-2)\bigr).
	\end{equation}
	We substitute~\eqref{eq:refined-known}, with $m=n-s-qb_k$, in the terms
	of~\eqref{eq:refined-rec}. This is what is done in the proofs of the two theorems, for the
	passage from $k-1$ to $k$; the only difference is that there, in Step~2, $s$ is replaced by
	one of its bounds $0$ and $b_k-a_k$, while here it is kept.

	\emph{Proof of \textup{(i)}.} Let $n\geq s$, and set $T\coloneqq n-s+M_{k-1}-(k-2)$, so that
	$m+M_{k-1}-(k-2)=T-qb_k$. Substituting the upper bound of~\eqref{eq:refined-known} in each
	term of~\eqref{eq:refined-rec}, and then adding the terms with
	$q>\floor*{\frac{n-s}{b_k}}$, which are non-negative (Lemma~\ref{lem:average}~(ii)), we get
	\begin{equation*}
		D_k(n)\leq\frac{1}{b_2\dotsm b_{k-1}}\sum_{q\geq0}\alpha_{k-1,k-2}(T-qb_k).
	\end{equation*}
	This is the sum $\sum_{q\geq0}\alpha_{k,k-1}(T-qb_{k+1})$ of Step~2 of the proof of
	Theorem~\ref{thm:main-upper}, with $k-1$ in place of $k$ and with the present $T$. From this
	point on the computation is that of the rest of Step~2 and of Steps~3 and~4 of that proof,
	which hold for every integer $T$; we do not repeat it. The result is
	\begin{equation*}
		D_k(n)
		\leq\frac{1}{b_2\dotsm b_k}\,\alpha_{k,k-1}\Bigl(T+\floor*{\frac{b_k}{2}}\Bigr)
		=\frac{1}{b_2\dotsm b_k}\,\alpha_{k,k-1}\bigl(n-s+M_k-(k-1)\bigr),
	\end{equation*}
	where the equality is $M_k=M_{k-1}+\floor*{\frac{b_k}{2}}+1$ (see~\eqref{eq:M-rec}). This
	is~\eqref{eq:claim-upper} with $n-s$ in place of $n$. Since $n-s\in\N$, the first part of
	the proof of Theorem~\ref{thm:main-upper}, with $n-s$ in place of $n$, shows that the
	right-hand side is the sum in~(i), and that it is at most
	$\frac{1}{b_2\dotsm b_k}\binom{n-s+M_k}{k-1}$.

	\emph{Proof of \textup{(ii)}.} Let $n\geq s+\Delta_{k-1}$, and set
	$T\coloneqq n-s-\Delta_{k-1}\in\N$, so that $m-\Delta_{k-1}=T-qb_k$. Substituting the lower
	bound of~\eqref{eq:refined-known} in each term of~\eqref{eq:refined-rec}, exchanging the two
	finite sums, and then adding the terms with $q>\floor*{\frac{n-s}{b_k}}$, which vanish
	because $T-qb_k<0$ for them, we get, by~\eqref{eq:S-def},
	\begin{equation*}
		D_k(n)
		\geq\frac{1}{b_2\dotsm b_{k-1}}\sum_{i=0}^{k-3}\bbinom[b-1]{k-3}{i}\sum_{q\geq0}f_{k-2-i}(T-qb_k)
		=\frac{1}{b_2\dotsm b_{k-1}}\sum_{i=0}^{k-3}\bbinom[b-1]{k-3}{i}S_{k-2-i,\,b_k}(T).
	\end{equation*}
	This is~\eqref{eq:lower-step2} with $k-1$ in place of $k$ and with the present $T$. From
	this point on the computation is that of Steps~3 and~4 of the proof of
	Theorem~\ref{thm:main-lower}, which hold for every integer $T$; we do not repeat it. The
	result is
	\begin{equation*}
		D_k(n)\geq\frac{1}{b_2\dotsm b_k}\sum_{i=0}^{k-2}\bbinom[b-1]{k-2}{i}f_{k-1-i}(n-s-\Delta_{k-1}).
	\end{equation*}
	Since $n-s-\Delta_{k-1}\geq0$, by Definition~\ref{dfn:trunc} we have
	$f_{k-1-i}(n-s-\Delta_{k-1})=\binom{n-s-\Delta_{k-1}+k-1-i}{k-1-i}$, and this is~(ii).
\end{proof}

\begin{rem}\label{rem:refined}
	Proposition~\ref{prop:refined} contains Theorems~\ref{thm:main-upper}
	and~\ref{thm:main-lower}, and is sharper. Recall that $0\leq s\leq b_k-a_k$.
	\begin{itemize}
		\item \emph{Upper bound.} If $n\geq s$, then $n-s\leq n$. The upper bounds of
		      Theorem~\ref{thm:main-upper} do not decrease with $n$, since they are
		      $\alpha_{k,k-1}$ and $f_{k-1}$ at the point $n+M_k-(k-1)$, up to the factor
		      $\frac{1}{b_2\dotsm b_k}$ (Lemmas~\ref{lem:average}~(ii)
		      and~\ref{lem:trunc}~(ii)). So~(i) implies Theorem~\ref{thm:main-upper}; if $n<s$,
		      the sum of Lemma~\ref{lem:recursion} is empty and $D_k(n)=0$.
		\item \emph{Lower bound.} By Lemma~\ref{lem:basic}~(iii) (for $k=2$, by
		      $\Delta_2=b_2-a_2$ and $\Delta_1=0$),
		      \begin{equation*}
			      s+\Delta_{k-1}\leq(b_k-a_k)+\Delta_{k-1}=\Delta_k .
		      \end{equation*}
		      Hence every $n\geq\Delta_k$ satisfies $n\geq s+\Delta_{k-1}$ and
		      $n-s-\Delta_{k-1}\geq n-\Delta_k$. Since the functions $f_{k-1-i}$ are
		      non-decreasing (Lemma~\ref{lem:trunc}~(ii)) and the coefficients
		      $\bbinom[b-1]{k-2}{i}$ are non-negative (Lemma~\ref{lem:sym-poly}~(ii)), (ii)
		      implies Theorem~\ref{thm:main-lower}. Moreover~(ii) gives $D_k(n)>0$ also for the
		      multiples $n$ of $d_k$, if any, with $s+\Delta_{k-1}\leq n<\Delta_k$, where
		      Theorem~\ref{thm:main-lower} says nothing.
	\end{itemize}
	The refinement concerns only the last step, for the reason explained in Remark~6.2 of the
	first paper: at the previous steps each term of the recursion has its own integer $s$, and
	the sums are no longer taken along one arithmetic progression.
\end{rem}

\subsection{Solutions of \texorpdfstring{$a_1x_1+\dotsb+a_kx_k\leq n$}{a1x1+...+akxk<=n}}
As in Section~6.2 of the first paper, for $k\in\N_{\geq1}$ and $n\in\N$ let
$\widehat{D}_k(n)$ be the number of $(x_1,\dotsc,x_k)\in\N^k$ with
$a_1x_1+\dotsb+a_kx_k\leq n$, and let
\begin{equation*}
	a'\coloneqq(a'_1,a'_2,a'_3,\dotsc)=(1,a_1,a_2,\dotsc),
	\qquad\text{that is}\qquad
	a'_1=1,\quad a'_i=a_{i-1}\ (i\geq2),
\end{equation*}
be the shifted sequence. We denote by $d'_i$, $b'_i$, $\Delta'_k$, $B'_k$, $M'_k$,
$\varepsilon'_i$, $E'_k$, $\bbinom[b'-1]{m}{i}$ and $\bbinom[-\varepsilon']{m}{i}$ the objects
built from $a'$ instead of $a$. Since $a'_1=1$, Section~\ref{s:a1}, applied to $a'$, gives
$d'_i=1$ and $b'_i=a'_i=a_{i-1}$ for every $i\geq2$, $\Delta'_k=0$ and $B'_k=-(k-1)$; hence
\begin{equation}\label{eq:shifted-constants}
	b'_2\dotsm b'_{k+1}=a_1\dotsm a_k,
	\qquad
	M'_{k+1}=a_1+\sum_{i=2}^{k}\Bigl(\floor*{\frac{a_i}{2}}+1\Bigr),
	\qquad
	\bbinom[b'-1]{k-1}{i}=e_i\Bigl(\frac{a_2-1}{2},\dotsc,\frac{a_k-1}{2}\Bigr).
\end{equation}
Moreover $\varepsilon'_i=1$ exactly when $a_{i-1}$ is even, so $E'_{k+1}$ is the number of even
numbers among $a_2,\dotsc,a_k$, and
$\bbinom[-\varepsilon']{k-1}{i}=\bigl(-\frac12\bigr)^{i}\binom{E'_{k+1}}{i}$
by~\eqref{eq:coeff-eps}.

\begin{cor}\label{cor:denum-ext}
	Let $k\in\N_{\geq1}$ and $n\in\N$, and set $N\coloneqq d_k\floor*{\frac{n}{d_k}}$. With the
	notation~\eqref{eq:shifted-constants},
	\begin{equation*}
		\frac{1}{a_1\dotsm a_k}\sum_{i=0}^{k-1}\bbinom[b'-1]{k-1}{i}\binom{N+k-i}{k-i}
		\leq\widehat{D}_k(n)
		\leq\frac{1}{a_1\dotsm a_k}\sum_{i=0}^{k-1}\bbinom[-\varepsilon']{k-1}{i}\binom{N+M'_{k+1}-i}{k-i}.
	\end{equation*}
	Moreover the lower bound is at least $\frac{1}{a_1\dotsm a_k}\binom{N+k}{k}$, and the upper
	bound is at most $\frac{1}{a_1\dotsm a_k}\binom{N+M'_{k+1}}{k}$.
\end{cor}
\begin{proof}
	\emph{Reduction to a denumerant.} By the first two steps of the proof of Corollary~6.4 of
	the first paper, $N$ is the largest multiple of $d_k$ not exceeding $n$,
	$\widehat{D}_k(n)=\widehat{D}_k(N)$, and $\widehat{D}_k(N)$ is the denumerant $D_{k+1}(N)$ of
	the sequence $a'$.

	\emph{The two bounds.} We have $k+1\geq2$, $d'_{k+1}=1$, and
	$N\geq0=B'_{k+1}+k$. So Theorem~\ref{thm:main-lower} and the
	bound~\eqref{eq:upper-sum} of Theorem~\ref{thm:main-upper}, applied to the sequence $a'$
	with $k+1$ in place of $k$, give
	\begin{equation*}
		\frac{1}{b'_2\dotsm b'_{k+1}}\sum_{i=0}^{k-1}\bbinom[b'-1]{k-1}{i}\binom{N-B'_{k+1}-i}{k-i}
		\leq\widehat{D}_k(N)
		\leq\frac{1}{b'_2\dotsm b'_{k+1}}\sum_{i=0}^{k-1}\bbinom[-\varepsilon']{k-1}{i}\binom{N+M'_{k+1}-i}{k-i},
	\end{equation*}
	which, with $-B'_{k+1}=k$ and $b'_2\dotsm b'_{k+1}=a_1\dotsm a_k$, are the two inequalities
	of the statement.

	\emph{The single binomial coefficients.} By Lemma~\ref{lem:sym-poly}, applied to $a'$, the
	term $i=0$ of the lower bound is $\binom{N+k}{k}$ and the other terms are non-negative. By
	the inequality~\eqref{eq:upper-single} of Theorem~\ref{thm:main-upper}, applied to the
	sequence $a'$ with $k+1$ in place of $k$, the upper bound is at most $\frac{1}{a_1\dotsm a_k}\binom{N+M'_{k+1}}{k}$.
\end{proof}

\begin{rem}\label{rem:denum-ext}
	The bounds of Corollary~\ref{cor:denum-ext} depend on $n$ only through $N$, as
	$\widehat{D}_k(n)$ does.
\end{rem}

\bigskip
\noindent\textbf{Acknowledgements.} The author used an AI assistant (Claude, Anthropic) to revise
the English wording and the organization of the text, and to check the formulas and the
numerical examples. All results and proofs are due to the author, who takes full responsibility
for the content of the paper.

\bibliographystyle{aomalpha}
\bibliography{denumerant_bounds_II}

\end{document}